\documentclass[journal,twoside,web]{ieeecolor}

\usepackage{amsmath,amsfonts,amssymb}
\usepackage{dsfont}
\usepackage{algorithm}
\usepackage{algpseudocode}
\usepackage{accents}
\usepackage{siunitx}
\usepackage{tikz}
\usepackage{pgfplots}
\pgfplotsset{compat=1.18}
\usepgfplotslibrary{groupplots}
\usepackage{cite}
\usepackage{enumerate}
\usepackage{mathtools}
\usepackage{generic}
\usepackage{graphicx}
\usepackage{textcomp}
	
\renewcommand{\baselinestretch}{0.95}

\newcommand*{\vm}[1]{\boldsymbol{#1}}
\newcommand*{\trans}{^{\top}}
\newcommand*{\dd}{\mathrm{d}}
\newcommand*{\nablax}{\nabla_{\vm x_i}}
\newcommand*{\nablaxx}{\nabla_{\vm x_i \vm x_i}^2}
\newcommand*{\nablaxj}{\nabla_{\vm x_j}}
\newcommand*{\neighs}{j \in \mathcal{N}_i}

\newcommand*{\agents}{i \in \mathcal{V}}
\newcommand*{\Ni}[1]{\vm{#1}_{\mathcal N_i}}
\newcommand*{\Nj}[1]{\vm{#1}_{\mathcal N_j}}
\newcommand*{\inds}{^{*}}
\newcommand*{\indqp}{^{q-1}}
\newcommand*{\indq}{^{q}}
\newcommand*{\st}{\operatorname{s.\!t.}}
\newcommand*{\smallspace}{\,\,}
\newcommand*{\ubar}[1]{\underaccent{\bar}{#1}}

\DeclareMathOperator*{\rank}{\mathrm{rank}}
\DeclareMathOperator*{\diag}{\mathrm{diag}}
\DeclareMathOperator*{\nullspace}{\mathrm{null}}

\newtheorem{theorem}{Theorem}
\newtheorem{lemma}{Lemma}
\newtheorem{assumption}{Assumption}
\newtheorem{corollary}{Corollary}
\newtheorem{remark}{Remark}

\newtheorem{example}{Example}

\newenvironment{contexample}{
	\addtocounter{example}{-1} \begin{example}[cont.]}{
\end{example}}

\def\BibTeX{{\rm B\kern-.05em{\sc i\kern-.025em b}\kern-.08em
T\kern-.1667em\lower.7ex\hbox{E}\kern-.125emX}}
\begin{document}
\title{Sensitivity-Based Distributed Programming for Non-Convex Optimization}
\author{Maximilian Pierer von Esch, Andreas V\"olz, and Knut Graichen, \IEEEmembership{Senior Member, IEEE}
\thanks{This work is funded by the Deutsche Forschungsgemeinschaft (DFG, German Research Foundation) under project no. 464391622.}
\thanks{The authors are with the Chair of Automatic Control, Friedrich-Alexander-Universität Erlangen-Nürnberg (FAU), Erlangen, Germany.
			Email: \{maximilian.v.pierer, andreas.voelz, knut.graichen\}@fau.de}}
	\maketitle
	%
	%%%%%%%%%%%%%%%%%%%%%%%%%%%%%%%%%%%%%%%%%%%%%%%%%%%%%%%%%%%%%%%%%%%%%%%%%%%%%%%%%%%%%%%%%%%%%%%%%%%%%%%%%%%%%%%
	%
	\begin{abstract}
		This paper presents a novel sensitivity-based distributed programming (SBDP) approach for non-convex, large-scale nonlinear programs (NLP). The algorithm relies on first-order sensitivities to cooperatively solve the central NLP in a distributed manner with only neighbor-to-neighbor communication and parallelizable local computations. The decoupling of the subsystems is based on primal decomposition. We derive sufficient local convergence conditions for non-convex problems. Furthermore, we consider the SBDP method in a distributed optimal control context and derive favorable convergence properties in this setting. We illustrate these theoretical findings and the performance of the proposed method with a comparison to state-of-the-art algorithms and simulations of various distributed optimization and control problems.
	\end{abstract}
	
	%
	%%%%%%%%%%%%%%%%%%%%%%%%%%%%%%%%%%%%%%%%%%%%%%%%%%%%%%%%%%%%%%%%%%%%%%%%%%%%%%%%%%%%%%%%%%%%%%%%%%%%%%%%%%%%%%%
	%
	\section{Introduction}
	\label{sec:intro}
	% Necessity of distributed optimization
	Many networked system applications, such as distributed model predictive control (DMPC) \cite{Christofides}, optimal power flow~\cite{Molzahn}, wireless sensor networks~\cite{Rabbat}, or machine learning~\cite{Boyd}, require solving large-scale optimization problems. Distributed optimization techniques address this issue by spreading computations across subsystems with limited information exchange, offering greater flexibility and computational scalability compared to centralized approaches.
	
	% Review of different optimization methods in particular sensitivity-based
	Consequently, a variety of algorithms have been proposed to decompose the central, large-scale optimization problem into lower-dimensional and parallelizable subproblems. One of the most popular methods is the alternating direction method of multipliers (ADMM) \cite{Boyd}, which shows promising results in various application areas. However, ADMM suffers from drawbacks, such as convergence guarantees for only certain special classes of non-convex problems~\cite{Hong,Houska} as well as requiring strong duality and an increased communication effort. In practice, many problems are inherently non-convex, e.g., in terms of DMPC with nonlinear systems \cite{Bestler,Burk} or power networks~\cite{Erseghe}, and the communication effort should be kept to a minimum, as it often presents the bottleneck in distributed optimization~\cite{Burk}.
	The augmented Lagrangian alternating direction inexact Newton (ALADIN) method is a well-known algorithm suitable for non-convex problems~\cite{Houska}. However, ALADIN employs a central coordination step that may be prohibitive regarding scalability. To alleviate this issue, bi-level variants have been proposed that distribute the coordination step with, e.g., ADMM~\cite{Engelmann}. In a similar spirit, algorithms based on distributing steps of classical centralized techniques like sequential quadratic programming (SQP)~\cite{Stomberg} or interior point methods~\cite{Engelmann2} have been developed. Another line of research regards establishing distributed versions of proximal minimization algorithms which are commonly used to address non-convexity in optimization \cite{Hong2,Jiang}. 

While most distributed optimization algorithms, such as ADMM, rely on a decomposition of the dual problem by introducing local copies of shared variables and enforcing consensus through dual updates, this approach can face challenges for non-convex nonlinear programs (NLPs). Convergence guarantees typically require convexity and strong duality~\cite{Boyd}. Furthermore, the size of the local subproblems is increased due to the auxiliary variables, which can lead to higher computational complexity and communication overhead. Primal decomposition methods, in contrast, iterate directly on the coupled decision variables. This avoids the need for strong duality and keeps subproblem dimensionality lower, which can reduce algorithmic complexity. However, ensuring the optimality of the distributed solution and feasibility of subproblems in the primal setting is generally non-trivial.

Besides Jacobi iterations for distributed, convex optimization \cite{Doan,Gross} or forward-backward splitting~\cite{Tang}, sensitivity-based approaches provide a promising option in this context. Rather than relying on dual variables, they exploit first-order sensitivities that quantify the effect of a neighbor’s decision on the local objectives and constraints. Prior work applied this idea to linear distributed MPC and state estimation \cite{Scheu,Scheu1,Alvarado,Schneider2,Schneider3}, while subsequent studies extend the concept to nonlinear systems \cite{Scheu2,Pierer,Pierer2}. However, convergence results for general nonlinear NLPs remain largely unexplored, and prior analysis either focus on linear-quadratic problems \cite{Scheu} or on a restricted class of nonlinear optimal control problems \cite{Pierer,Pierer2}. This highlights that the sensitivity-based approach is still largely unexplored for general non-convex NLPs with coupled constraints.

	%Contribution and structure
	Therefore, this work focuses on sensitivity-based distributed programming (SBDP) for non-convex problems defined over networks with a generic topology. The approach relies on first-order sensitivities, which quantify the interaction between subsystems to cooperatively solve the centralized NLP and can be computed efficiently in a distributed manner. In this aspect, SBDP constitutes a first-order method in the sense that no explicit Hessians or higher-order derivatives are needed. The approach is characterized by solving lower-dimensional NLPs at the subsystem level and requiring only neighbor-to-neighbor communication. Furthermore, we derive the SBDP method both for a general class of NLPs and for NLPs with a specific neighbor-affine structure~\cite{Burk}. The latter is particularly advantageous, as it allows formulating the SBDP scheme with only one communication step per iteration, compared to two in the general case. Since one of the main intended areas of application is nonlinear DMPC, we demonstrate the usage and investigate the convergence of SBDP for optimal control. 
	%Contribution
	The contribution of this work can be summarized as follows:
	\begin{itemize}
	\item  We extend the concept of SBDP from convex, linear-quadratic problems \cite{Scheu} to a class of NLPs with non-convex objectives and (coupled) non-convex constraints.
	\item  We rigorously investigate the local convergence properties and show that SBDP is at least linearly convergent under appropriate conditions. Additionally, we derive criteria, depending on the problem structure and interaction between subsystems, for SBDP to converge quadratically.
	\item We explicitly consider SBDP in a distributed optimal control context. In particular, we explore the effect of the prediction horizon length on its convergence behavior.
	\item We validate the theoretical results with several numerical simulations of relevant problems and demonstrate the efficacy of the proposed approach through the distributed optimal control of coupled inverted pendulums.
	\end{itemize}
	
	%Structure
	The paper is structured in the following way: Section~\ref{sec:ProblemStatement} introduces the problem formulation, for which the SBDP approach is derived in Section \ref{sec:dist_solution}. Local convergence of SBDP is analyzed in Section \ref{sec:alg_analysis}. In Section \ref{sec:DMPC}, the application of the approach to nonlinear optimal control problems is discussed in more detail. The numerical evaluation is presented in Section \ref{sec:num_eval}, before the paper is summarized in Section \ref{sec:conclusion}.
	
	% Notation 
	Notation: Given a vector $\vm v \in \mathbb{R}^n $ or matrix $\vm M \in \mathbb{R}^{n\times m}$ and an integer $i$, the notations $[\vm v]_i$ and $[\vm M]_i$ refer to the $i$-th component and $i$-th row of $\vm v$ and $\vm M$, respectively, while $\|\cdot\|$ denotes the respective Euclidean or induced spectral norm. Similarly, given an index set $\mathcal{S} \subset \mathbb{N}$, $[\vm M]_{\mathcal{S}}$ denotes the matrix consisting of the rows $[\vm M]_i$, $i\in \mathcal{S}$. The stacking of individual column vectors $\vm v_i \in \mathbb{R}^{n_i}$ , $i \in \mathcal{S}$ from an ordered index set $\mathcal{S} \subset \mathbb{N}$ is stated as $[\vm v_i]_{i \in \mathcal{S}}$. The maximum and minimum eigenvalues of	a square matrix $\vm M$ are denoted by $\lambda_{\min}(\cdot)$ and $\lambda_{\max}(\cdot)$. The set of integers in the range from $0$ to $N$ is denoted as $\mathbb{I}_{[0,N]}$. An $r$-neighborhood of a point $\vm v_0 \in \mathbb{R}^v$ is defined as $\mathcal{B}_r(\vm v_0):=\{ \vm v \in \mathbb{R}\,| \|\vm v - \vm v_0\| < r\}$. For the Jacobian of a vector-valued function $\vm f(\vm x) : \mathbb{R}^n \rightarrow \mathbb{R}^m$, we write $\nabla \vm f(\vm x) = [\nabla f_1(\vm x), \dots, \nabla f_m(\vm x)]\trans \in \mathbb{R}^{m\times n}$, whereby $\nabla f_i$ is the gradient of the $i$-th component. 
	%
	%%%%%%%%%%%%%%%%%%%%%%%%%%%%%%%%%%%%%%%%%%%%%%%%%%%%%%%%%%%%%%%%%%%%%%%%%%%%%%%%%%%%%%%%%%%%%%%%%%%%%%%%%%%%%%%%%%%%%%%%%%%%%%%%%%%%%%%%%%%%%%%%%%%%%%%%%%%%%%%%%%%%%%%%%%%%%%%%%%%%%%
	%
	\section{Problem Statement}
	\label{sec:ProblemStatement}
	We consider NLPs whose structure is described by a connected, undirected graph $\mathcal{G}=(\mathcal{V},\mathcal{E})$ where the set of vertices $\mathcal{V}=\{1,\dots,M\}$ denotes a number of subsystems, referred to as agents, and the set of edges $\mathcal{E}\subset \mathcal{V}\times \mathcal{V}$ defines couplings between agents. The agents aim at cooperatively solving the central NLP
	\begin{subequations}\label{eq:central_NLP}
		\begin{align}
			\min_{\vm x_1,\dots,\vm x_M} &\quad  \sum_{i\in\mathcal V} f_i(\vm x_i, \Ni{x}) \label{eq:central_costFunction}\\
			 ~\st \quad&\quad \vm  g_i( \vm x_i,\Ni{x}) = \vm 0\,, \quad \agents \label{eq:central_equality} \\
			&\quad \vm h_i(\vm x_i, \Ni{x})\leq \vm 0\,, \quad \agents \label{eq:central_inequality}
		\end{align}
	\end{subequations}
	with the local decision variables $\vm x_i \in \mathbb{R}^{n_i}$ which are optimized by each agent $\agents$. The coupling to the neighbors is given by the stacked notation $\Ni{x}:=[\vm x_j]_{\neighs}$ with $\mathcal{N}_i:=\{j\in \mathcal{V}: (i,j)\in \mathcal{E},\, i\neq j\}$ denoting the set of neighbors for every agent $\agents$. Furthermore, each agent minimizes the objective functions $f_i:\mathbb{R}^{n_i}\times\mathbb{R}^{\sum_{\neighs} n_j} \rightarrow \mathbb{R}$ subject to (coupled) equality constraints $\vm g_i:\mathbb{R}^{n_i}\times\mathbb{R}^{\sum_{\neighs} n_j} \rightarrow \mathbb R^{n_{gi}}$ and (coupled) inequality constraints $\vm h_i:\mathbb{R}^{n_i}\times \mathbb{R}^{\sum_{\neighs} n_j} \rightarrow \mathbb R^{n_{hi}}$, where for each $\agents$, $n_i$ is the number of local decision variables, $n_{gi}$ is the number of local equality constraints and $n_{hi}$ is the number of local inequality constraints. All functions in NLP \eqref{eq:central_NLP} are assumed to be at least three times continuously differentiable. 
	We define the central Lagrangian of~\eqref{eq:central_NLP} as
	\begin{align}
		L(\vm x, \vm \lambda, \vm \mu) = \sum_{\agents}L_i(\vm x_i, \vm \lambda_i, \vm \mu_i, \Ni{x}) \label{eq:central_Lagrangian}
	\end{align}
	and the local Lagrangians $L_i = L_i(\vm x_i, \vm \lambda_i, \vm \mu_i,\Ni{x})$ as
	\begin{align}
		L_i:= f_i(\vm x_i, \Ni{x}) + \vm \lambda_i\trans \vm  g_i( \vm x_i,\Ni{x}) + \vm \mu_i\trans  \vm h_i(\vm x_i, \Ni{x}) \label{eq:local_Lagrangian}
	\end{align}
	for every $\agents$. The quantities  $\vm \lambda_i \in \mathbb{R}^{n_{gi}}$ and $\vm \mu_i \in \mathbb{R}^{n_{hi}}$ in \eqref{eq:local_Lagrangian} represent the Lagrange multipliers associated with the respective constraints \eqref{eq:central_equality} and \eqref{eq:central_inequality} of the central NLP. The centralized notations $\vm x =[\vm x_i]_{\agents} \in \mathbb R^{n}$, $ [\vm \lambda_i]_{\agents} \in \mathbb{R}^{n_g}$, and $\vm \mu = [\vm \mu_i]_{\agents} \in \mathbb{R}^{n_h}$ summarize all local variables and are concatenated as  $ \vm p:=[\vm x_i\trans,\, \vm \lambda_i\trans,\, \vm \mu_i\trans ]\trans_{\agents} \in \mathbb{R}^p$ with $p=n+n_g+n_h$, where $n$ is the number of global decision variables, $n_g$ is the number of global equality constraints, and $n_h$ is the number of global inequality constraints. In case of large-scale systems, NLP \eqref{eq:central_NLP} represents a high dimensional optimization problem. Therefore, the remainder of this paper is concerned with the distributed solution of NLP~\eqref{eq:central_NLP}. 
	%
	%%%%%%%%%%%%%%%%%%%%%%%%%%%%%%%%%%%%%%%%%%%%%%%%%%%%%%%%%%%%%%%%%%%%%%%%%%%%%%%%%%%%%%%%%%%%%%%%%%%%%%%%%%%%%%%%%%%%%%%%%%%%%%%%%%%%%%%%%%%%%%%%%%%%%%%%%%%%%%%%%%%%%%%%%%%%%%%%%%%%%%
	%
	\section{Distributed sensitivity-based solution}
	\label{sec:dist_solution}
	The main idea of the proposed SBDP approach for the distributed solution of the central NLP \eqref{eq:central_NLP} involves augmenting the individual cost functions $f_i(\vm x_i, \Ni{x})$ in \eqref{eq:central_costFunction} of each agent $\agents$ with the first-order sensitivities of their neighbors $\neighs$ and resolving the couplings in the cost function~\eqref{eq:central_costFunction} and constraints \eqref{eq:central_equality}-\eqref{eq:central_inequality} via a primal decomposition approach. First, we address the solution of the general NLP \eqref{eq:central_NLP} via SBDP, before we discuss a special case in which objective and constraint functions exhibit specific structural properties.
	\subsection{SBDP for a general structure}
Each agent $\agents$ constructs a local, decoupled NLP which is solved in each iteration $q=1,2,\dots$ of SBDP 
	\begin{subequations}\label{eq:local_NLP}
		\begin{align}
			(\vm x_i\indq,\, \vm \lambda_i\indq,\, \vm \mu_i\indq) = \arg\min_{\vm x_i} \quad & \phi_i(\vm x_i)
			\\  ~\st \quad& \vm  g_i( \vm x_i,\Ni{x}\indqp) = \vm 0 \label{eq:local_equality} \\
			& \vm h_i(\vm x_i, \Ni{x}\indqp)\leq \vm 0\,, \label{eq:local_inequality}
		\end{align}
	\end{subequations}
	with the local, modified cost function $\phi_i:\mathbb{R}^{n_i} \rightarrow \mathbb{R}$
	\begin{equation} \label{eq:local_costFunction}
		\phi_i(\vm x_i)\!:=\!f_i(\vm x_i, \Ni{x}\indqp) + \sum_{\neighs}\big(\nablax  L_j\indqp \big)\trans(\vm x_i - \vm x_i\indqp)\,,
	\end{equation}
	which implicitly depends on the iterates at $q-1$. The quantities $(\vm \lambda_i\indq, \vm \mu_i\indq)$ denote the dual variables associated with the constraints \eqref{eq:local_equality} and \eqref{eq:local_inequality}, respectively, and are assumed to be available from the solution of~\eqref{eq:local_NLP}. Otherwise, they may be recovered from the corresponding KKT conditions of~\eqref{eq:local_NLP}.
The first-order sensitivity term in \eqref{eq:local_costFunction} takes the influence of the local decision variables on the neighboring agents objective function into account.  It is defined as the directional derivative of the local Lagrangian $L_j(\cdot)$, $\neighs$, given by the gradient of the previous iteration $  \nablax L_j\indqp:= \nablax L_j(\vm x_j\indqp, \vm \lambda_j\indqp, \vm \mu_j\indqp, \Nj{x}\indqp)$, in direction of $\vm x_i - \vm x_i\indqp$. The couplings are resolved in a primal decomposition fashion by using the decision variables $\Ni{x}\indqp$ of the previous iteration which need to be communicated. The gradient $\nablax L_j(\cdot)$ is computed from \eqref{eq:local_Lagrangian} as
	\begin{align}\label{eq:gradient_withoutneighboraffine}
		& \nablax L_j(\vm x_j, \vm \lambda_j, \vm \mu_j, \Nj{x}) =  \nablax f_j(\vm x_j, \Nj{x}) \nonumber                        \\
		& \phantom{=}+ \nablax \vm g_j(\vm x_j, \Nj{x})\trans \vm \lambda_j + \nablax \vm h_j(\vm x_j, \Nj{x})\trans \vm \mu_j\,,
	\end{align}
where $\vm \lambda_j$ and $\vm \mu_j$ are the Lagrange multipliers belonging to the respective neighbor, $\neighs$. Note that \eqref{eq:gradient_withoutneighboraffine} may potentially involve decision variables of second-order neighbors, i.e., $ \mathcal{N}_i^2 := (\bigcup_{\neighs} \mathcal{N}_j)\setminus (\mathcal{N}_i \cup \{i\})$, which are not directly available in a neighbor-to-neighbor communication network. This is addressed by agent $i$ computing the mirroring gradient $\nablaxj L_i(\cdot)$ instead and communicating it to the neighbors.
The decoupling of the central NLP \eqref{eq:central_NLP} into the local NLPs~\eqref{eq:local_NLP} is exploited by solving the individual problems in parallel at the agent level, see Algorithm~\ref{alg:SENSI_without_neighboraffinity}. This requires the agents to communicate bi-directionally with their neighbors. Step 1 involves the evaluation of the partial derivatives $\nablaxj L_i\indqp$ at the iterates of the previous iteration $q-1$, which are subsequently communicated in Step~2 to the neighbors such that the sensitivity term can be appended to the cost function \eqref{eq:local_costFunction}. Afterward, NLP \eqref{eq:local_NLP} is solved locally and in parallel by each agent before the new decision variable $\vm x_i\indq$ is communicated to the respective neighbors. A possible stopping criterion is the difference between two iterates, i.e.,
	\begin{equation} \label{eq:stopping_criterion}
	\|\vm p\indq - \vm p\indqp\|\leq \epsilon
	\end{equation} 
with tolerance $\epsilon>0$. If $\|\cdot\|_\infty$ is chosen, then \eqref{eq:stopping_criterion} can be evaluated in a distributed fashion and only convergence flags must be communicated. Algorithm \ref{alg:SENSI_without_neighboraffinity} relies exclusively on local computations and two neighbor-to-neighbor communication steps per SBDP iteration, with at most $\sum_{\agents}2n_i|\mathcal{N}_i|$ floats exchanged system-wide, resulting in a distributed scheme that scales well with the number of agents \cite{Pierer3}.
\begin{remark}
The sensitivity term in \eqref{eq:gradient_withoutneighboraffine} may be approximated using numerical differentiation schemes when analytical gradients are not available. While this is a viable option, SBDP may converge only to an approximate central solution, whose optimality depends on the accuracy of the approximation.
\end{remark}
	%
	% General SBDP 
	\begin{algorithm}[tb]\small
		\caption{General SBDP for each agent $\agents$}
		\begin{algorithmic}[1]
			\Statex Initialize $\vm x_i^0$, $\vm \lambda_i^0$, $\vm \mu_i^0$, send $\vm x_i^0$ to neighbors $\neighs$,  set $q\leftarrow 1$
			\State Compute $\nablaxj  L_i\indqp$ via \eqref{eq:gradient_withoutneighboraffine} for all $\neighs$
			\State Send  $\nablaxj  L_i\indqp$ to all $\neighs$
			\State Compute $\vm x_i^{q}$, $\vm \lambda_i^{q}$ and $\vm \mu_i^{q}$ by solving NLP \eqref{eq:local_NLP} to local optimality
			\State Send $\vm x_i^{q}$ to neighbors $\neighs$
			\State Stop if a suitable convergence criterion, e.g. \eqref{eq:stopping_criterion}, is met. Otherwise, return to line~$1$ with $q \leftarrow q+1$.
		\end{algorithmic}\label{alg:SENSI_without_neighboraffinity}
	\end{algorithm}
	\subsection{SBDP for a neighbor-affine structure}
	An especially favorable case arises if the objective and constraint functions exhibit a certain structure, that is, if they can be expressed in a so-called neighbor-affine form \cite{Burk}. This means that the coupling to the neighbors in NLP \eqref{eq:central_NLP} enters additively via nonlinear coupling functions which depend on exactly one other decision variable $\vm x_j,\, \neighs$ such that the objective and constraint functions \eqref{eq:central_costFunction}-\eqref{eq:central_inequality} are of the form
\begin{align}\label{eq:central_NLP_neighboraffine}
		\vm k_i(\vm x_i, \Ni{x}) &= \vm k_{ii}(\vm x_i)+ \sum_{\neighs} \vm k_{ij}(\vm x_i, \vm x_j)
\end{align}
with appropriately dimensioned local functions $\vm k_{ii}(\cdot)$ and pairwise coupling functions $\vm k_{ij}(\cdot,\cdot)$.
 Evaluating the gradient $\nablax L_j(\cdot)$ for this case leads to a simplified version of~\eqref{eq:gradient_withoutneighboraffine}, i.e.,
	\begin{align}\label{eq:gradient_neighboraffine}
		& \nablax L_j(\vm x_j, \vm \lambda_j, \vm \mu_j, \Nj{x}) =  \nablax f_{ji}(\vm x_j, \vm x_i) \nonumber                            \\
		& \phantom{=}+ \nablax \vm g_{ji}(\vm x_j, \vm x_i)\trans \vm \lambda_j + \nablax \vm h_{ji}(\vm x_j, \vm x_i)\trans \vm \mu_j\,.
	\end{align}
	The key takeaway is that the agents are now able to compute the sensitivities locally since the right-hand side of \eqref{eq:gradient_neighboraffine} only depends on local and neighboring variables and not on variables of  second-order neighbors, $\mathcal{N}_i^2$, as in \eqref{eq:gradient_withoutneighboraffine}. This property is leveraged to reduce the number of communication steps by sending the Lagrange multipliers together with the decision variables in Step 2 of Algorithm~\ref{alg:SENSI_with_neighboraffinity}, which summarizes the modified SBDP method for neighbor-affine functions. This results in only one neighbor-to-neighbor communication step in which a maximum number of $\sum_{\agents}(n_i+n_{gi} +n_{hi}) |\mathcal{N}_i|$ floats needs to be sent system-wide per iteration. While the number of communicated floats is not necessarily lower than in Algorithm~\ref{alg:SENSI_without_neighboraffinity}, the communication overhead, e.g., (de)constructing data or waiting times, resulting from an extra communication step is nullified. For decoupled constraints, the multipliers can be omitted in Step~2, as they do not affect the gradient \eqref{eq:gradient_neighboraffine}, thereby reducing the amount of exchanged data.
	\begin{remark}
		Neighbor-affine structures frequently occur in various problems in the area of DMPC, e.g., water networks \cite{Burk,Alvarado} or robotic formations \cite{Rosenfelder}. In addition, many optimal power flow problems are inherently neighbor-affine~\cite{Erseghe}.
	\end{remark}
	%
	% SBDP for neighbor-affine 
	\begin{algorithm}[tb]\small
		\caption{Neighbor-affine SBDP for each agent $\agents$}
		\begin{algorithmic}[1]
			\Statex Initialize and send $\vm x_i^0$, $\vm \lambda_i^0$, $\vm \mu_i^0$ to $\neighs$, set $q\leftarrow 1$
			\State Compute $\vm x_i^{q}$, $\vm \lambda_i^{q}$ and $\vm \mu_i^{q}$ by solving NLP \eqref{eq:local_NLP} to local optimality
			\State Send $\vm x_i^{q}$, $ \vm \lambda_{i}\indq$ and $\vm \mu_{i}^{q}$ to neighbors $\neighs$
			\State Stop if a suitable convergence criterion, e.g. \eqref{eq:stopping_criterion}, is met. Otherwise, return to line~$1$ with $q \leftarrow q+1$. 
		\end{algorithmic}\label{alg:SENSI_with_neighboraffinity}
	\end{algorithm}
	%  
	% Example
	\begin{example} 
	\textit{As an example consider the non-convex NLP
		\begin{subequations}\label{eq:example_constrained_NLP}
			\begin{align}
				\min_{x_1,\,x_2} \quad & x_1^2(x_1^2-2)+x_2^2(x_2^2-2)+x_1^2x_2^2\label{eq:example_constrained_NLP_cost_function}
				\\  ~\st \quad& 3x_1 -x_2 - 2 = 0\,, \label{eq:example_constrained_NLP_equality}
			\end{align}
		\end{subequations}
		with the two decision variables $x_1,\, x_2 \in \mathbb{R}$ which are to be optimized by two agents $i\in \mathcal{V}=\{1,2\}$. If we partition the problem as $f_{ii} = x_i^2(x_i^2-2)$, $f_{ij} =0.5x_i^2x_j^2$, $g_1= g_{11} +g_{12}$ with $g_{11}=3x_1-2$ and $g_{12}=-x_2$, it is in neighbor-affine form \eqref{eq:central_NLP_neighboraffine}. By computing the gradients according to \eqref{eq:gradient_neighboraffine}, i.e.,
		\begin{align}
		\nabla_{x_1}L_2 = x_1x_2^2\,,\quad \nabla_{x_2}L_1 = x_2x_1^2 - \lambda_1\,,
		\end{align}
		each agent $\agents$ can construct the modified NLPs \eqref{eq:local_NLP} which can be solved in parallel according to Algorithm \ref{alg:SENSI_with_neighboraffinity}.}
	\end{example}
	%%%%%%%%%%%%%%%%%%%%%%%%%%%%%%%%%%%%%%%%%%%%%%%%%%%%%%%%%%%%%%%%%%%%%%%%%%%%%%%%%%%%%%%%%%%%%%%%%%%%%%%%%%%%%%%%%%%%%%%%%%%%%%%%%%%%%%%%%%%%%%%%%%%%%%%%%%%%%%%%%%%%%%%%%%%%%%%%%%%%%%
	%
	\section{Algorithmic Analysis}
	\label{sec:alg_analysis}
	The algorithmic analysis of Algorithms \ref{alg:SENSI_without_neighboraffinity} and \ref{alg:SENSI_with_neighboraffinity} can be performed jointly as they only differ in the computation of the sensitivities. The analysis is organized into two sections, of which the first states the central and local Karush-Kuhn-Tucker (KKT) conditions together with
	necessary regularity assumptions. The second section investigates the convergence towards a central KKT point and contains the proof of convergence of Algorithms \ref{alg:SENSI_without_neighboraffinity} and \ref{alg:SENSI_with_neighboraffinity}.
	\subsection{Central and local optimality conditions}
	The KKT conditions of the central NLP~\eqref{eq:central_NLP} are given as 
	\begin{subequations}\label{eq:centralKKT}
		\begin{alignat}{2}
			 \vm 0 &= \nabla_{\vm x_i} L(\vm x, \vm \lambda, \vm \mu) \,,  &&\agents \label{eq:centralKKT_gradient} \\
			\vm 0 &= \vm g_i(\vm x_i, \Ni{x})\,, \hspace{1mm}&&\agents\\
			\vm 0 & = \vm U_i \vm h_i(\vm x_i, \Ni{x}) ,\,\vm h_i(\vm x_i, \Ni{x})\leq\vm  0,\, \vm \mu_i \geq \vm 0,\,\,&& \agents
		\end{alignat}
	\end{subequations}
	with the diagonal matrix $\vm U_i:= \diag(\mu_{1,i},\dots,\mu_{n_{hi},i})$ for each $\mu_{k,i}$ in $\vm \mu_i$, $k \in\mathbb{N}_{[1,n_{hi}]}$. 
	Let $\vm p\inds:= [(\vm x_i\inds)\trans,\, (\vm \lambda_i\inds)\trans,\, (\vm \mu_i\inds)\trans ]_{\agents}\trans$ be a KKT point of~\eqref{eq:central_NLP} which satisfies~\eqref{eq:centralKKT}. The set $\mathcal{A}(\vm x)$ of active inequality constraints is given by 
	\begin{align}
			\mathcal{A}(\vm x) := \{k \in\mathbb{N}_{[1,n_h]}\, |\,  [\vm h(\vm x)]_k = 0 \} 
		\end{align}
with the notation $\vm g(\vm x):=[\vm g_i(\vm x_i, \Ni{x})]_{\agents}$ and $\vm h(\vm x):=[\vm h_i(\vm x_i, \Ni{x})]_{\agents}$ summarizing all constraints.
	The equalities of the central KKT conditions~\eqref{eq:centralKKT} are equivalently written as
	\begin{equation} \label{eq:centralKKTstacked}
		\vm F(\vm p) := \begin{bmatrix}
			\nabla_{\vm x_i} L(\vm x, \vm \lambda,\vm \mu ) \\
			\vm g_i(\vm x_i,\Ni{x})                         \\
			\vm U_i \vm h_i(\vm x_i,\Ni{x})
		\end{bmatrix}_{\agents} = \vm 0
	\end{equation}
	with $\vm F: \mathbb{R}^p \rightarrow \mathbb{R}^p$.
The following standard regularity assumption regarding the KKT point of problem~\eqref{eq:central_NLP} is made.
	\begin{assumption}\label{ass:regularity}
		There exists a KKT point $\vm p\inds$ of \eqref{eq:central_NLP} which
		\begin{enumerate}[i)]
			\item satisfies $ [\vm \mu\inds]_k> 0$ if $[\vm h(\vm x\inds)]_k= 0$, $k\in \mathbb{N}_{[1,n_h]}$, i.e., strict complementary slackness (SCS) holds,
			\item satisfies the second-order sufficient conditions (SOSC), i.e., $\vm s\trans \nabla_{\vm x \vm x}^2L(\vm x \inds, \vm \lambda\inds, \vm \mu\inds) \vm s>0$ for all $\vm s \neq \vm 0 $ with $\nabla_{\vm x} \vm g(\vm x\inds) \vm s= \vm 0$ and $[\nabla_{\vm x} \vm h(\vm x \inds )]_{\mathcal{A}(\vm x\inds)} \vm s= \vm 0$,
			\item satisfies the linear constraint qualification (LICQ), i.e., the column vectors of $[\nabla_{\vm x} \vm g(\vm x\inds)\trans, [\nabla_{\vm x} \vm h(\vm x\inds)]_{\mathcal{A}(\vm x\inds)}\trans]$ are linearly independent.
		\end{enumerate}
	\end{assumption}
	Regarding the modified local NLPs \eqref{eq:local_NLP}, the corresponding KKT conditions are established as
	\begin{subequations} \label{eq:localKKT}
		\begin{align}
			 \vm 0 &= \nablax L_i(\vm x_i, \vm \lambda_i, \vm \mu_i,\Ni{x}\indqp) + \sum_{\neighs} \nablax L_j\indqp                    \label{eq:localKKT_gradient} \\
		 \vm 0 &= \vm g_i(\vm x_i, \Ni{x}\indqp)                                                                                                   \\
		 \vm 0 &= \vm U_i \vm h_i(\vm x_i, \Ni{x}\indqp),\,	\vm h_i(\vm x_i,\Ni{x}\indqp)\leq\vm  0,\, \vm \mu_i \geq \vm 0 \label{eq:localKKT_inequality}
		\end{align}
	\end{subequations}
for every $\agents$ in each iteration $q$. We define the KKT point $\vm p\indq := [\vm p_i\indq]_{\agents}$ with $\vm p_i:= [\vm x_i\trans, \vm \lambda_i\trans, \vm \mu_i\trans]\trans$ as the solution of all KKT systems \eqref{eq:localKKT} in each iteration $q$. The set 
	\begin{align}
			\mathcal{A}_i(\vm x_i, \Ni{x}):= \{k \in \mathbb{N}_{[1,n_{hi}]}\, |\, [\vm h_i(\vm x_i,\Ni{x})]_k = 0  \} 
	\end{align}
	denotes the active constraints of the modified NLPs~\eqref{eq:local_NLP}. Similar to \eqref{eq:centralKKTstacked}, we write
	\begin{align} \label{eq:localKKTstacked}  
		 &\vm {\hat F}(\vm p, \vm p\indqp) :=  \\      
		  &\begin{bmatrix}
			\nablax L_i(\vm x_i, \vm \lambda_i, \vm \mu_i, \Ni{x}\indqp)+ \sum_{\neighs} \nablax L_j\indqp \\
			\vm g_i(\vm x_i, \Ni{x}\indqp)                                                                 \\
			\vm U_i \vm h_i(\vm x_i, \Ni{x}\indqp) \nonumber
		\end{bmatrix}_{\agents} = \vm 0
	\end{align}
	to summarize the equalities of the local KKT conditions \eqref{eq:localKKT} with the function $\vm {\hat F}: \mathbb{R}^p \times \mathbb{R}^p \rightarrow \mathbb{R}^p $. This reveals that the modified NLPs~\eqref{eq:local_NLP} are parametric in $[\vm p_j\indqp]_{\neighs \cup \{i\}}$ or equivalently in $\vm p\indqp$ from a global viewpoint. Consequently, we introduce the (implicit) mapping $ \vm \Phi: \mathbb{R}^p \rightarrow \mathbb{R}^p $ as
	\begin{equation}\label{eq:mapping_Phi}
	\vm p\indq = \vm \Phi(\vm p\indqp)\,, \quad \vm {\hat F}(\vm \Phi(\vm p\indqp), \vm p\indqp)= \vm 0\,.
	\end{equation}
	In general, the mapping $\vm \Phi(\cdot)$ cannot be written in closed form and might even be a set-valued solution associated with the KKT system \eqref{eq:localKKT}.
	However, under certain regularity conditions, $\vm \Phi(\cdot)$ is locally single-valued and continuously differentiable in a neighborhood of the centralized KKT point $\vm p\inds$ such that \eqref{eq:mapping_Phi} holds. Then, the local single-valuedness allows us to analyze the SBDP iteration as a fixed-point recursion on $\vm{\Phi}(\cdot)$ in a neighborhood of $\vm p\inds$. The required regularity conditions are stated below and can be seen as a compatibility condition between problems \eqref{eq:central_NLP} and \eqref{eq:local_NLP}.
	\begin{assumption} \label{ass:compatibility}
		For all $\agents$, the central KKT point $\vm p\inds$ of NLP \eqref{eq:central_NLP} satisfies
		\begin{enumerate}[i)]
			\item the SOSC for \eqref{eq:local_NLP}, i.e., $\vm s_i\trans \nablaxx L_i(\vm x_i\inds, \vm \lambda_i\inds, \vm \mu_i\inds, \Ni{x}\inds) \vm s_i>0$ for all $\vm s_i \neq \vm 0 $ with $\nablax\vm g_i(\vm x_i \inds, \Ni{x}\inds) \vm s_i= \vm 0$ and $[\nablax  \vm h_i(\vm x_i \inds, \Ni{x}\inds )]_{\mathcal{A}_i(\vm x_i\inds, \Ni{x}\inds)} \vm s_i= \vm 0$,
			\item the LICQ for \eqref{eq:local_NLP}, i.e., the column vectors of the matrix $[\nablax \vm g_i(\vm x_i\inds,\vm x_{\mathcal{N}_i}\inds)\trans, [\nablax \vm h_i(\vm x_i\inds,\vm x_{\mathcal{N}_i}\inds)]_{\mathcal{A}_i(\vm x_i\inds,\vm x_{\mathcal{N}_i}\inds)}\trans]$ are linearly independent.
		\end{enumerate}
	\end{assumption}
	Assumption \ref{ass:compatibility} requires a proper assignment of the objective functions and constraints to the modified NLPs \eqref{eq:local_NLP} such that $\vm p\inds$ satisfies the SOSC and LICQ at the agent level.  In this sense, it presents a slightly stricter version of Assumption \ref{ass:regularity} and generalizes \cite[Ass. 3 and 4]{Scheu} to the nonlinear setting. If instead a weaker constraint qualification than LICQ is used, the local Lagrange multipliers may not be unique and hence $\vm \Phi(\cdot)$ might be set-valued \cite[Thm. 2]{Wachsmuth}.
	\begin{remark}
If Statement i) in Assumption \ref{ass:compatibility} is not satisfied, a remedy is to add the quadratic term $\rho/2 \| \vm x_i - \vm x_i\indqp\|^2$ with penalty parameter $\rho\geq0$ to \eqref{eq:local_costFunction}. By selecting $\rho\geq \max\{ |\lambda_{\min} (\nablaxx L_i(\vm x_i\inds, \vm \lambda_i\inds, \vm \mu_i\inds, \Ni{x}\inds))|:\agents\}$ with the minimum eigenvalue $\lambda_{\mathrm{min}}(\cdot)$, i) is guaranteed to hold without affecting the distributed minimizer.  
		\end{remark}
		
	\begin{remark}
			If Statement ii) in Assumption \ref{ass:compatibility} is not satisfied, it may be possible to repartition the central NLP~\eqref{eq:central_NLP}
			such that ii) holds. Moreover, ii) implies that the number of active constraints at $\vm x\inds$ must
			be equal or less than $n$. If $n_{gi} + |\mathcal{A}_i(\vm x_i\inds, \Ni{x}\inds)|>n_i$ or the local LICQ is not satisfied, the introduction of slack variables or local copies with corresponding consistency constraints may serve as a remedy to ensure the feasibility of the local NLPs~\eqref{eq:local_NLP}.
		\end{remark}
	
	\subsection{Local convergence}
	\label{subsec:localConvergence}
	The subsequent analysis will revolve around showing that the sequence $\{\vm p\indq\}$ generated by repeatedly solving the KKT system \eqref{eq:localKKT} converges towards a (local) central solution $\vm p\inds$ of the KKT system \eqref{eq:centralKKT}.
To associate the limit of this sequence with the KKT point of \eqref{eq:central_NLP}, it is necessary to show that $\vm\Phi(\vm p\inds)$, i.e., the concatenated KKT point of all modified problems \eqref{eq:local_NLP} for $ \vm p \indqp = \vm p\inds$, is also a KKT point of the central NLP \eqref{eq:central_NLP}. This relation is established in the next lemma.
\begin{lemma}\label{lem:optConditions}
Let Assumptions \ref{ass:regularity} and \ref{ass:compatibility} hold. Then, any point $\vm {p}$ which satisfies  $ \vm p = \vm \Phi(\vm {p})$ together with the corresponding inequalities in \eqref{eq:localKKT_inequality} also satisfies~\eqref{eq:centralKKT}. In particular, $ \vm p\inds= \vm \Phi(\vm p\inds)$ is a locally unique primal-dual solution of \eqref{eq:localKKT}.
	\end{lemma}
	\begin{proof}
		The assertion directly follows from the equivalence of the  KKT conditions \eqref{eq:centralKKT} and \eqref{eq:localKKT}, where the following relation $\nabla_{\vm x_i} L(\vm x, \vm \lambda, \vm \mu) = \nablax L_i(\vm x_i, \vm \lambda_i, \vm \mu_i, \Ni{x}) + \sum_{\neighs} \nablax L_j(\vm x_j, \vm \lambda_j, \vm \mu_j, \Nj{x})$ is utilized to show the equality of the right-hand sides of \eqref{eq:centralKKT_gradient} and \eqref{eq:localKKT_gradient} for all $\agents$. Thus, the relation
			$\vm F(\vm p) = \vm {\hat F}(\vm \Phi(\vm p), \vm p)$
together with the equivalence of the inequalities $\vm h(\vm x) = [\vm h_i(\vm x_i, \Ni{x})]_{\agents}$ and $\vm \mu = [\vm \mu_i ]_{\agents} \geq \vm 0$ imply that $\vm \Phi(\vm p)$ is also a KKT point of~\eqref{eq:central_NLP}. Inserting $\vm p\indqp = \vm p\inds$ into \eqref{eq:localKKT}, shows that $\vm p\inds $ is a solution for these equations under the validity of \eqref{eq:centralKKT}. Local uniqueness follows from Assumptions \ref{ass:regularity} and \ref{ass:compatibility}, cf. \cite[Thm. 3.2.2]{Fiacco}.
	\end{proof}
Lemma \ref{lem:optConditions} shows that the local KKT system \eqref{eq:localKKT} has a locally unique solution at $ \vm p\inds = \vm \Phi(\vm p\inds)$ that satisfies the central KKT conditions \eqref{eq:centralKKT}. This ensures the optimality the distributed solution.
Assumptions \ref{ass:regularity} and \ref{ass:compatibility}, the result of Lemma~\ref{lem:optConditions}, and the interpretation of the modified NLPs~\eqref{eq:local_NLP} as a perturbed version of the central NLP \eqref{eq:central_NLP} can be used to invoke the basic sensitivity theorem \cite[Thm. 3.2.2]{Fiacco} which is adapted for the present case and summarized in Lemma~\ref{lem:solvability}.

\begin{lemma}\label{lem:solvability}
		Under Assumptions \ref{ass:regularity} and \ref{ass:compatibility}, there exists a constant $r_1>0$ such that for $\vm p \in \mathcal{B}_{r_1}(\vm p\inds)$ it holds that
		\begin{enumerate}[i)]
			\item the $C^1$ mapping $\vm \Phi(\cdot)$ is locally single-valued with $\vm \Phi(\vm p)$ satisfying $\vm{\hat F}(\vm \Phi(\vm p), \vm p)=\vm 0$ and the inequalities in \eqref{eq:localKKT_inequality},
			\item at any $\vm \Phi(\vm p)$ the set of active inequalities of NLP \eqref{eq:central_NLP} and of the local NLPs \eqref{eq:local_NLP} are equivalent, the LICQ holds for the NLPs~\eqref{eq:local_NLP}, SCS is preserved and the SOSC holds.
		\end{enumerate}
\end{lemma}
	\begin{proof}
		Lemma \ref{lem:optConditions} shows that for $\vm p\indqp = \vm p\inds$ the local solutions $\vm p\indq$ of NLP~\eqref{eq:local_NLP} satisfy the central KKT conditions~\eqref{eq:centralKKT}. Furthermore, Assumption \ref{ass:compatibility} implies that the second-order sufficiency conditions and the LICQ hold for each NLP~\eqref{eq:local_NLP}, $\agents$ at $\vm p\inds$. Strict complementarity for each NLP \eqref{eq:local_NLP} at $\vm p\inds$ for all inequality constraints follows from the equality of $\vm h(\vm x \inds) = [\vm h_i(\vm x_i\inds, \Ni{x}\inds)]_{\agents}$ and Statement i) in Assumption~\ref{ass:regularity}. Thus, with all involved functions being at least three times continuously differentiable, the conditions of the basic sensitivity theorem \cite[Thm. 3.2.2]{Fiacco} are satisfied, allowing its application to each NLP~\eqref{eq:local_NLP} yielding the statements above.
	\end{proof}
In essence, Lemma \ref{lem:solvability} ensures that the mapping $\vm \Phi(\cdot)$ is differentiable on $\mathcal{B}_{r_1}(\vm p\inds)$ and that the regularity assumptions of the central KKT point $\vm p\inds$ carry over to the KKT point $\vm p\indq = \vm \Phi(\vm p \indqp)$ of the modified NLPs \eqref{eq:local_NLP} such that their feasibility is ensured in every iteration, provided $\vm p\indqp \in \mathcal{B}_{r_1}(\vm p\inds)$. The following lemma exploits the differentiability of $ \vm \Phi(\cdot)$ to explicitly compute an expression of its Jacobian.
	\begin{lemma}\label{lem:jacobian}
		Suppose Assumptions \ref{ass:regularity} and \ref{ass:compatibility} hold. Then, the Jacobian  $\vm J(\vm p):=\nabla \vm \Phi(\vm p)$ of the mapping $\vm \Phi(\vm p)$, implicitly defined in \eqref{eq:mapping_Phi},  for $\vm p \in \mathcal{B}_{r_1}(\vm p\inds)$ is given by
		\begin{equation} \label{eq:Jacobian}
			\vm J(\vm p)= -\vm M(\vm p)^{-1} \vm N(\vm p)\,,
		\end{equation}
		where the explicit expressions of $\vm M(\vm p) \in \mathbb{R}^{p\times p}$ and $\vm N(\vm p) \in \mathbb{R}^{p\times p} $ are stated in Appendix \ref{app:Jacobian}.
	\end{lemma}
	\begin{proof}
		The expression \eqref{eq:Jacobian} is obtained by computing the total derivative of the function $\vm{\hat F}(\vm q, \vm p) = \vm{\hat F}(\vm \Phi(\vm p),\vm p)$, as defined in \eqref{eq:localKKTstacked}, w.r.t. $\vm p\in \mathcal{B}_{r_1}(\vm p\inds)$. Since $\vm {\hat F}(\vm \Phi(\vm p), \vm p) = \vm 0$ must hold for all $\vm p \in \mathcal{B}_{r_1}(\vm p\inds)$, it must also hold for the total derivative that $\frac{\mathrm{d}\vm {\hat F}}{\mathrm{d}\vm p} =\vm 0$. This results in the expression
		\begin{align} \label{eq:Jacobian_derivation}
		\nabla_{\vm q} \vm {\hat F}(\vm \Phi(\vm p), \vm p)\, \nabla \vm \Phi(\vm p) + \nabla_{\vm p} \vm {\hat F}(\vm \Phi(\vm p), \vm p)  =  \vm 0\,.
		\end{align}
		With the abbreviation $\vm M(\vm p):=\nabla_{\vm q} \vm {\hat F}(\vm \Phi(\vm p), \vm p)$ and $ \vm N(\vm p):= \nabla_{\vm p}\vm {\hat F}(\vm \Phi(\vm p), \vm p)$ one directly obtains \eqref{eq:Jacobian}. The regularity of $\vm M(\vm p)$ also follows from \cite[Thm. 3.2.2]{Fiacco}.
	\end{proof}

\begin{remark}
The structure of $\vm J(\vm p)$ exposes the necessity of the local SOSC and LICQ in Assumption \ref{ass:compatibility} which ensures that the matrix $\vm M(\vm p)$ is regular. If, e.g., the LICQ does not hold, the mapping $\vm \Phi(\cdot)$ might not be differentiable in $\mathcal{B}_{r_1}(\vm p\inds)$, impeding the algorithmic analysis.
	\end{remark}
	\begin{contexample}
	\textit{
	Returning to the example NLP \eqref{eq:example_constrained_NLP}, the Jacobian \eqref{eq:Jacobian} can be found by at first computing $\vm M(\vm p)$ and $\vm N(\vm p)$ according to \eqref{eq:M_i} and \eqref{eq:N} in Appendix \ref{app:Jacobian}  
	\begin{align} \nonumber
	\vm M(\vm p) \!=\!\! \begin{bsmallmatrix}\!
	12x_1^2 + x_2^2-4 & 2 & 0  \\
	2 & 0 & 0 \\
	0&0&12x_2^2 + x_1^2-4
	\end{bsmallmatrix},\,
	\vm N(\vm p) \! = \! \begin{bsmallmatrix}
		x_2^2 & 0 & 4x_1x_2  \\
		0 & 0 & -1 \\
		4x_1x_2&-1&x_1^2
	\end{bsmallmatrix}
	\end{align} 
	and then evaluating the expression \eqref{eq:Jacobian} to arrive at 
	\begin{align}\label{eq:Jacobian_exampleNLP}
	\vm J(\vm p)\!=\! -\begin{bsmallmatrix}
		0 & 0 & -\frac{1}{3}  \\
		\frac{1}{3}x_2^2 &0 & \frac{1}{3}(4x_1^2 + \frac{1}{3}x_2^2+4x_1x_2-\frac{4}{3})\\
		\frac{4x_1x_2}{12x_1^2 +x_2^2-4}&\frac{-1}{12x_1^2 +x_2^2-4}&	\frac{x_1^2}{12x_1^2 +x_2^2-4}
	\end{bsmallmatrix}
	\end{align}
	which, in this case, is independent of the multiplier $\lambda_1$ since the coupled inequality constraint function \eqref{eq:example_constrained_NLP_equality} is linear.
	}
	\end{contexample}
	Although the mapping $\vm \Phi(\cdot)$ is generally not available in closed-form, its Jacobian $\vm J(\vm p)$ can be used to construct an approximation of the next iterate from the previous iterate in a neighborhood of $\vm p\inds$. This fact is used to establish local convergence of SBDP as stated in the following theorem.
	\begin{theorem} \label{th:linconv}
		Let Assumptions \ref{ass:regularity} and \ref{ass:compatibility} hold. Then, there exists a constant $r>0$ such that if $\vm p^0 \in \mathcal{B}_r(\vm p\inds)$ and $\|\vm J(\vm p^*)\|<1$, the sequence $\{\vm p\indq\}$ generated by Algorithm~\ref{alg:SENSI_without_neighboraffinity} or \ref{alg:SENSI_with_neighboraffinity} converges to a central KKT point $\vm p\inds$ of NLP~\eqref{eq:central_NLP}. Moreover, for all $q=1,2,\dots$, the following holds:
		\begin{enumerate}[i)]
			\item The order of convergence is at least linear, i.e., there exists a constant $C\in(0,1)$ such that $\|\vm p\indq - \vm p\inds\| \leq C\|\vm p\indqp - \vm p\inds\|$ holds,
			\item If $\|\vm J(\vm p\inds)\|=0$, then the order of convergence is quadratic, i.e., there exists a constant $\tilde C<\infty$ such that $\|\vm p\indq - \vm p\inds\| \leq \tilde C\|\vm p\indqp - \vm p\inds\|^2$ holds.
		\end{enumerate}
	\end{theorem}
	\begin{proof}
		See Appendix \ref{app:Th1}.
	\end{proof}
	Theorem \ref{th:linconv} is constructive as it provides a worst case estimate of the convergence rate and the size of  $\mathcal{B}_r(\vm p\inds)$ such that in practice the performance will usually be better.
	\begin{remark}\label{remark:estimate_convergenceRadius}
		A (conservative) estimate of the convergence radius \( r \) can be obtained by finding the largest \( r \) such that \( \frac{L}{2} \|\vm p^0 - \vm p\inds\| + \|\vm J(\vm p\inds)\| < 1 \) for all \( \vm p^0 \in \mathcal{B}_r(\vm p\inds) \), which follows directly from \eqref{eq:ineqJ}. Here, \( r \) must be chosen such that the (local) Lipschitz constant \( L \) of the Jacobian is finite (\( L < \infty \)) on \( \mathcal{B}_r(\vm p\inds) \), and \( r < r_1 \) to ensure the applicability of Lemma \ref{lem:solvability}. 
	\end{remark}
	The condition $\|\vm J(\vm p^*)\|<1$ is a requirement on the ratio between the change of the KKT conditions~\eqref{eq:localKKT} caused by the local decision variables and that induced by neighboring ones and can be interpreted as a generalized diagonal dominance condition on the NLP~\eqref{eq:central_NLP}. If this condition is not met, the steps taken by the algorithm in the local optimal direction, i.e., solving the NLPs~\eqref{eq:local_NLP}, might be too large to result in a global contraction. This requirement is consistent with results from the linear-quadratic problems \cite{Scheu,Schneider2,Schneider3} and shows that for highly coupled systems the algorithm might diverge since then the external change in the local KKT conditions~\eqref{eq:localKKT} could exceed the limit for contraction. As suggested in \cite{Schneider2}, it might be possible to re-partition the problem such that $\|\vm J(\vm p\inds )\|<1$ is satisfied, if the partitioning of NLP \eqref{eq:central_NLP} is treated as a design parameter and not as fixed.
	Furthermore, Theorem \ref{th:linconv} states that the order of convergence is at least linear, which, however, can be strengthened by imposing additional requirements on the Jacobian $\vm J(\vm p)$. In particular, Statement ii) of Theorem \ref{th:linconv} requires that $\|\vm J(\vm p\inds)\|=0$. Although this is certainly not the case for all problem classes, an example with quadratic convergence is part of the numerical simulations and can be explicitly shown for certain problems, as discussed in the following Corollary.
	\begin{corollary}
		Let Assumptions \ref{ass:regularity} and \ref{ass:compatibility} hold. Furthermore, suppose that NLP \eqref{eq:central_NLP} is of neighbor-affine form \eqref{eq:central_NLP_neighboraffine} with objective functions and constraints that satisfy $\nabla^2_{\vm x_i\vm x_j}f_{ij} = \nabla^2_{\vm x_j\vm x_j}f_{ij} = \vm 0$, $ \nabla_{\vm x_j} \vm g_{ij} = \nabla^2_{\vm x_i\vm x_j} \vm g_{ij} = \nabla^2_{\vm x_j\vm x_j} \vm g_{ij} = \vm 0$ and $\nabla_{\vm x_j} \vm h_{ij} = \nabla^2_{\vm x_i\vm x_j} \vm h_{ij} = \nabla^2_{\vm x_j\vm x_j} \vm  h_{ij} = \vm 0$ at the minimizer $\vm x\inds$. Then, the order of convergence of Algorithm \ref{alg:SENSI_with_neighboraffinity} is quadratic. 
	\end{corollary}
	\begin{proof}
		The claim directly follows from $\vm N(\vm p\inds) = \vm 0$ according to \eqref{eq:N} in Appendix \ref{app:Jacobian} which implies that $\vm J(\vm p\inds) = \vm 0$ holds such that Statement ii) of Theorem~\ref{th:linconv} is applicable.
	\end{proof}
	\begin{remark}
		For decoupled local constraints, the iterates $\vm x\indq$ are feasible in each iteration $q$ w.r.t. the central constraints \eqref{eq:central_equality}-\eqref{eq:central_inequality}. This is particularly interesting for suboptimal DMPC, where feasibility can be used to establish stability~\cite{Scokaert}. This is a major advantage over dual approaches, which usually achieve primal feasibility only in the limit.
	\end{remark}
	\begin{remark}
		Although only first-order information (i.e., gradients of the local Lagrangians) is exchanged between neighboring agents, each agent internally solves a nonlinear program and thereby exploits local higher-order information in its update. In this sense SBDP combines first-order communication with second-order processing at the agent level, which explains the local quadratic convergence for certain problem structures. This fundamentally differs from pure first-order schemes such as distributed gradient descent, where both the exchanged and the locally processed information are genuinely first-order.
	\end{remark}
	\begin{contexample}
	\textit{
		On the affine constraint subspace $\mathbb{G}=\{(x_1,x_2)\in \mathbb{R}^2\, :\, g_1(x_1,x_2)=0 \}$, the cost function~\eqref{eq:example_constrained_NLP_cost_function} possesses a local minimum at $\vm x\inds = [0.35,-0.94]\trans$ with multiplier $\lambda_1\inds = 0.20$. The regularity Assumptions \ref{ass:regularity} and \ref{ass:compatibility} are clearly satisfied at $\vm x\inds$. Evaluating the Jacobian \eqref{eq:Jacobian_exampleNLP} at the local minimum, results in $\|\vm J(\vm p\inds)\|=0.76$ which shows that Algorithm \ref{alg:SENSI_with_neighboraffinity} converges at a linear order of convergence with a rate of $C= 0.76$ when $\|\vm p\indq - \vm p\inds\|$ is small. According to Remark \ref{remark:estimate_convergenceRadius}, a (conservative) estimate of the convergence radius is computed as $r\approx 0.35$, where the (local) Lipschitz constant $L=1.43$ of \eqref{eq:Jacobian_exampleNLP} in $\mathcal{B}_{0.35}(\vm p\inds)$ was computed numerically. 
	}
	\end{contexample}
	%
	%%%%%%%%%%%%%%%%%%%%%%%%%%%%%%%%%%%%%%%%%%%%%%%%%%%%%%%%%%%%%%%%%%%%%%%%%%%%%%%%%%%%%
	%%%%%%%%%%%%%%%%%%%%%%%%%%%%%%%%%%%%%%%%%%%%%%%%%%%%%%%%%%%%%%%%%%%%%%%%%%%%%%%%%%%%%
	%
	\section{Application to distributed optimal control}
	\label{sec:DMPC}
	The proposed SBDP method has proven itself to be a powerful tool in the context of cooperative DMPC \cite{Alvarado,Pierer,Scheu} and its counterpart, distributed moving horizon estimation \cite{Schneider2,Schneider3}, where it solves the central optimal control problem in a distributed fashion. Therefore, we analyze the SBDP method within the context of nonlinear optimal control and examine its convergence behavior from a control perspective.
	
	\subsection{Central optimal control problem}
	\label{subsec:OCP}
	In this setting, each agent $\agents$ is additionally associated with the nonlinear, continuous-time system dynamics
	\begin{equation}\label{eq:agent_dynamics}
		\vm{\dot {\ubar x}}_i(t) = \vm {f}_i(\vm {\ubar x}_i(t), \vm u_i(t), \Ni{\ubar x}(t))\,,\smallspace \forall t>0
	\end{equation}
	with the local states $\vm {\ubar x}_i(t)\in \mathbb{R}^{n_{xi}}$, controls $\vm u_i(t)\in \mathbb{R}^{n_{ui}}$, initial conditions $\vm {\ubar x}_i(0)= \vm{ \ubar x}_{i,0}$ and functions $ \vm {f}_i:\mathbb{R}^{n_{xi}}\times \mathbb{R}^{n_{ui}} \times \mathbb{R}^{\sum_{\neighs} n_{xj}}\rightarrow \mathbb{R}^{n_{xi}}$. Additionally, we assume the existence and uniqueness of the respective solution $\vm {\ubar x}_i(\cdot;\vm{ \ubar x}_{i,0},\vm u_i(\cdot), \Ni{\ubar x}(\cdot))$ for every $\agents$. The underscore notation highlights the difference between local optimization vector $\vm x_i$ and local state $\vm{\ubar{x}}_i(t)$. Similar to~\eqref{eq:central_NLP}, the coupling to the neighbors $\neighs$ in \eqref{eq:agent_dynamics} is given by the states $\vm{ \ubar x}_j(t)$. Then, the agents aim at cooperatively solving the central optimal control problem (OCP) in a distributed fashion
	\begin{subequations}\label{eq:DMPC_cont_OCP}
		\begin{align}
		\min_{\vm u(\cdot) } \quad & \sum_{i\in\mathcal V} V_i(\vm {\ubar x}_i(T)) + \int_0^T\! l_i( \vm {\ubar x}_i(t), \vm u_i(t), \Ni{ \ubar x}(t))\, \dd t\label{eq:DMPC_cont_OCP_costFunction}
		\\  ~\st \quad& \vm{\dot {\ubar x}}_i(t) = \vm {f}_i(\vm {\ubar x}_i(t), \vm u_i(t), \Ni{\ubar x}(t))\,, \quad \agents \\
		&\vm {\ubar x}_i(0)= \vm{ \ubar x}_{i,0}\,, \quad \agents \\
		& \vm h_i(\vm {\ubar x}_i(t), \vm u_i(t))\leq \vm 0,\, \smallspace \vm h_i(\vm {\ubar x}_i(T))\leq \vm 0\,,\smallspace \agents\,, \label{eq:DMPC_cont_OCP_inequality}
		\end{align}
	\end{subequations}
    where the inequalities \eqref{eq:DMPC_cont_OCP_inequality} are to be understood pointwise-in-time and where $\vm u(t) = [\vm u_i(t)]_{\agents} $ denotes the stacked central control vector. The cost functional \eqref{eq:DMPC_cont_OCP_costFunction}  with the prediction horizon $T>0$ consists of the terminal cost $V_i:\mathbb{R}^{n_{xi}}\rightarrow \mathbb{R}$ and the integral cost $ l_i:\mathbb{R}^{n_{xi}}\times \mathbb{R}^{n_{ui}} \times \mathbb{R}^{\sum_{\neighs}n_{xj}}\rightarrow \mathbb{R}$.
To transform the problem into the form of~\eqref{eq:central_NLP}, 	we discretize the time interval $[0,T]$ equidistantly with the discrete time points $t_k = k \Delta t$, $k \in \mathbb{N}_0$ into $N \in \mathbb{N}$ subintervals with length $\Delta t = T/N>0$ and parameterize the control on the subintervals as piecewise constant, i.e., $\vm u_i(t) = \vm u_i^k$, $\forall t \in [t_k,t_{k+1})$. The state representation at $t_k$ is given as $\vm x_i^{k}=\vm {\ubar x}_i(t_k)$. To preserve the continuous-time coupling structure of \eqref{eq:agent_dynamics}, the states of the neighbors are kept constant on the discretization interval, i.e., $\Ni{x}^k = \Ni{ \ubar x}(t) $, $\forall t \in [t_k,t_{k+1})$. This simplification also preserves the neighbor-affine structure of the continuous-time problem if the dynamics \eqref{eq:agent_dynamics} and integral costs \eqref{eq:DMPC_cont_OCP_costFunction} are in neighbor-affine form \eqref{eq:central_NLP_neighboraffine}. A discretized version of the continuous-time OCP \eqref{eq:DMPC_cont_OCP} reads
	\begin{subequations}\label{eq:DMPC_OCP}
		\begin{align}
			\min_{\vm{ \hat x},\, \vm{\hat  u} } \quad & \sum_{i\in\mathcal V} V_i(\vm x_i^N) +  \Delta t \sum_{k=0}^{N-1} l_i(\vm x_i^k, \vm u_i^k, \Ni{x}^k) \label{eq:DMPC_OCP_costFunction}
			\\  ~\st \quad& \vm x_i^{k+1} =\vm x_i^k \!+\!\Delta t\vm f_i^d(\vm x_i^k, \vm u_i^k,\Ni{x}^k),\, k \in \mathbb{I}_{[0,N-1]} \label{eq:DMPC_OCP_equality} \\
			& \vm x_i^0 = \vm{ \ubar x}_{i,0}  \label{eq:DMPC_OCP_init} \\
			& \vm h_i(\vm x_i^k, \vm u_i^k)\leq \vm 0,\, k\in \mathbb{I}_{[0,N-1]}, \smallspace \vm h_i(\vm x_i^N)\leq \vm 0\,, \label{eq:DMPC_OCP_inequality}
		\end{align}
	\end{subequations}
	where the (in)equalities \eqref{eq:DMPC_OCP_equality}-\eqref{eq:DMPC_OCP_inequality} are to be understood for each $\agents$. The notation $\vm{\hat x}_i := [\vm x_i^k]_{k \in \mathbb{I}_{[0,N]}} \in \mathbb{R}^{(N+1)n_{x_i}}$, $\vm{\hat u}_i := [\vm u_i^k]_{k \in \mathbb{I}_{[0,N-1]}}\in \mathbb{R}^{Nn_{u_i}}$, $\vm{\hat x} := [\vm {\hat x}_i]_{\agents}$ and $\vm{\hat u} := [\vm {\hat u}_i]_{\agents}$ summarizes the decision variables.
	Typically, the discretization function $\vm {f}_i^d(\cdot)$  is derived from a numerical integration scheme. Furthermore, we consider only coupling in state variables in \eqref{eq:agent_dynamics} and \eqref{eq:DMPC_cont_OCP_costFunction} with subsystem-exclusive control variables to maintain a reasonable complexity and notation. The following considerations are equally applicable if input couplings are present in \eqref{eq:agent_dynamics} and \eqref{eq:DMPC_cont_OCP_costFunction} or if an implicit integration scheme is used to obtain \eqref{eq:DMPC_OCP_costFunction} and \eqref{eq:DMPC_OCP_equality}. The discretized OCP \eqref{eq:DMPC_OCP} is of the form \eqref{eq:central_NLP} and can be solved with Algorithms \ref{alg:SENSI_without_neighboraffinity} or \ref{alg:SENSI_with_neighboraffinity}.

	\subsection{Convergence analysis}
	\label{subsec:conv_analysis_OCP}
	In the following, we investigate the convergence of Algorithms \ref{alg:SENSI_without_neighboraffinity} and \ref{alg:SENSI_with_neighboraffinity} for the NLP \eqref{eq:DMPC_OCP} as a special case of \eqref{eq:central_NLP}. In particular, the following theorem establishes that a carefully chosen prediction horizon $T$ can guarantee convergence regardless of the coupling structure.
	\begin{theorem}\label{th:OCP_case}
		Suppose Assumptions \ref{ass:regularity} and \ref{ass:compatibility} hold for all $T\in \mathbb{T}$ on some compact and non-empty interval $T \in \mathbb{T}$. Then, for a given $N$, there exist a constant $r>0$ and a maximum prediction horizon $T_{\mathrm{max}}>0$ such that if $\vm p^0 \in \mathcal{B}_r(\vm p\inds)$, $ T<T_{\mathrm{max}}$ and $T_{\mathrm{max}} \in \mathrm{int}(\mathbb{T})$ the sequence $\{\vm p\indq\}$ generated by Algorithm~\ref{alg:SENSI_without_neighboraffinity} or \ref{alg:SENSI_with_neighboraffinity} converges to a KKT point $\vm p\inds$ of OCP~\eqref{eq:DMPC_OCP}.
	\end{theorem}
	\begin{proof}
		See Appendix \ref{app:Th2}. 
	\end{proof}
	Theorem \ref{th:OCP_case} shows that an upper bound on the prediction horizon $T_{\mathrm{max}}$ can always be found such that for $T<T_{\mathrm{max}}$ convergence of the SBDP method is guaranteed as long as the KKT point $\vm p\inds$ corresponding to $T_{\mathrm{max}}$ does not violate Assumptions \ref{ass:regularity} and \ref{ass:compatibility}. 
	Theorem \ref{th:OCP_case} can be strengthened if it is known a-priori that Assumptions \ref{ass:regularity} and \ref{ass:compatibility} hold on the complete optimization interval $T \in (0,\infty)$. Then, it can be guaranteed that an upper bound $T_{\mathrm{max}}$ always exists. 

	Considering SBDP in the context of DMPC provides the opportunity to utilize the prediction horizon not only as a design parameter for the DMPC scheme but also as an additional SBDP tuning parameter to enforce convergence. Depending on the MPC formulation, it is well known that the prediction horizon influences the stability of the DMPC scheme \cite{Jadbabaie,Limon,Graichen}. In light of Theorem~\ref{th:OCP_case}, it is a natural choice to select an MPC formulation without additional terminal constraints for which reducing $T$ might jeopardize the feasibility of~\eqref{eq:DMPC_OCP}. For a certain domain of attraction, stability can instead be ensured through a suitable terminal cost acting as a control Lyapunov function \cite{Jadbabaie,Limon,Graichen}. Then, the choice of $T$ represents a trade-off between convergence and the size of the domain of attraction. However, it does not (directly) influence the stability of the DMPC scheme. The maximum allowable horizon length can often be significantly increased by damping the iterates \cite{Pierer2}.
	In case of simple problems, the bound $T_{\mathrm{max}}$ can be explicitly computed as discussed in the next example. 
	\begin{example} 
		\textit{
			Consider two agents $i\in\mathcal{V}=\{1,\,2\}$ with the integrator dynamics $ \dot x_i(t) = u_i(t)$, integral cost functions $l_{i} = qx_i^2(t) + ru_i^2(t) + wx_i(t)x_j(t)$, weights $q,\,r>0$ and $w\in \mathbb{R}$ as well as zero terminal cost. Coupled cost functions of this kind frequently appear for example in the formation control of mobile robots \cite{Rosenfelder}. A forward Euler approximation of the agent dynamics \eqref{eq:agent_dynamics}, i.e., $\vm f_i^d(\cdot)= \vm f(\cdot)$, with $N=2$ yields the the discrete-time OCP
			\begin{subequations}\label{eq:example_OCP}
				\begin{align}
					\min_{\vm{\hat x},\, \vm{\hat u}} \quad & \Delta t\sum_{i\in\mathcal V}   l_i(x_i^0, u_i^0, x_{j}^0) + l_i(x_i^1, u_i^1, x_{j}^1) \label{eq:example_OCP_costFunction}
					\\  ~\st \quad& x_i^{k+1} = x_i^{k} + \Delta t u_i^k\,, \quad   k\in\{0,1\}\label{eq:example_OCP_equality} \\
					& x_i^0 = \ubar x_{i,0}\,,
				\end{align}
			\end{subequations}
			for some initial condition $\ubar x_{i,0} \in \mathbb{R}$.
			For this linear-quadratic example, the theoretical bound on the horizon $T_{\mathrm{max}}$ can be determined in closed form by computing the maximum eigenvalue $\lambda_{\mathrm{max}}(\vm J(\vm p))$ of the Jacobian  and imposing the convergence condition $|\lambda_{\mathrm{max}}(\vm J(\vm p))|<1$, i.e.,
			$
			|\lambda_{\mathrm{max}}(\vm J(\vm p))| = \big|\frac{2w\Delta t^2}{r+q \Delta  t^2} \big |<1\,.
			$
Thus, only if $2|w|>q$, we obtain
\begin{equation} \label{eq:T_max_example_OCP}
			T<T_{\mathrm{max}}= \sqrt{\frac{4r}{2|w|-q}}\,,
\end{equation} 
			as an upper bound on the horizon, otherwise the SBDP method is always convergent since $\lambda_{\mathrm{max}}(\vm J(\vm p))<1$ holds independently of $T$. Note that the condition \eqref{eq:T_max_example_OCP} is necessary and sufficient as \eqref{eq:example_OCP} is an equality constrained QP.
		}
	\end{example}	
	\section{Numerical evaluation}
	\label{sec:num_eval}
	Three numerical examples are utilized to analyze the performance of the SBDP method. The first concerns a simple non-convex minimization problem to illustrate the convergence properties derived in Theorem \ref{th:linconv} and compares the scheme to ADMM and ALADIN. The second addresses a more complex problem of the distributed optimal control of multiple coupled inverted pendulums, while the third investigates the role of the prediction horizon on the convergence behavior of SBDP. 
	\subsection{Convergence properties}
	\label{subsec:conv_properties}
\begin{figure}[tb]
		\centering
		\includegraphics{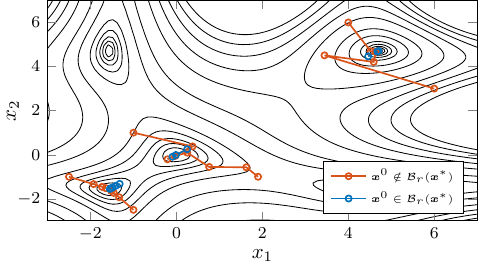}
		\vspace{-4mm}
		\caption{Contour lines of NLP \eqref{eq:example_NLP} and iteration profiles of Algorithm \ref{alg:SENSI_with_neighboraffinity}. The initializations of the blue iterates lie within $\mathcal{B}_r(\vm x\inds)$, while those in red do not. However, the SBDP method still converges for the initializations outside $\mathcal{B}_r(\vm x\inds)$ as the practical convergence radius is typically much larger.}
		\label{fig:conv_steps}
	\end{figure}

\begin{figure}[tb]
		\centering
		\includegraphics{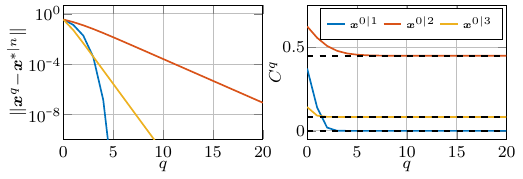}
		\vspace{-4mm}
		\caption{Error progression for different initial guesses (left) and empirical convergence rate $C^q$ (right) of Algorithm \ref{alg:SENSI_with_neighboraffinity} for NLP \eqref{eq:example_NLP}. The dashed lines show the norm of the Jacobian, $\|J(\vm x^{*|n})\|$, at the respective local minima~$n$.}
		\label{fig:conv_curves}
	\end{figure}

\begin{figure}[tb]
		\centering
		\includegraphics{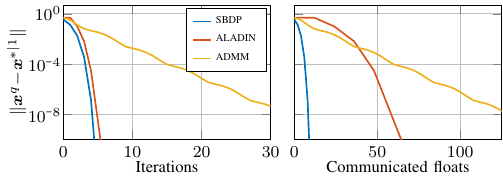}
		\vspace{-4mm}
		\caption{Comparison of SBDP, ADMM and ALADIN in terms of error progression to $ \vm x^{*|1}$ vs. iterations (left) and communicated floats (right).}
		\label{fig:compare_curves}
	\end{figure}

	For the evaluation of the convergence properties, we consider the following non-convex, unconstrained NLP
	\begin{equation}\label{eq:example_NLP}
		\min_{x_1,\, x_2} f(\vm x)= x_1^2 + x_2^2\sin x_1 + x_2^2 + x_1^2 \sin x_2
	\end{equation}
	with the two decision variables $x_1, x_2\in \mathbb{R}$, $\vm x = [x_1,x_2]\trans$, which are to be optimized by the agents $i\in \mathcal{V}=\{1,2\}$, respectively. Problem \eqref{eq:example_NLP} is in neighbor-affine form, such that we partition it according to \eqref{eq:central_NLP_neighboraffine}, where $f_{ii}(x_i) = x_i^2$ and $f_{ij}(x_i,x_j) = x_j^2\sin x_i$, leading to $L_i(x_i,x_j) = f_{ii}(x_i) +f_{ij}(x_i,x_j)$, $\agents$, $\neighs$. We compute \eqref{eq:gradient_neighboraffine} as $\nablax L_j = 2x_i\sin x_j$ and construct the NLPs \eqref{eq:local_NLP}, which simplify to
	\begin{align} \label{eq:example_localNLP}
		\min_{x_i} \smallspace x_i^2 +(x_j^{q-1})^{2}\sin x_i +2x_i\indqp \sin x_j\indqp (x_i - x_i\indqp)
	\end{align}
	for each agent and apply Algorithm \ref{alg:SENSI_with_neighboraffinity}, where we solve the local problems \eqref{eq:example_localNLP} via bisectioning. To illustrate the results of Theorem \ref{th:linconv}, we inspect the convergence behavior in the neighborhood of the following three local minima $\vm x^{*|n} \in\{[0,0]\trans, -[\frac{\pi}{2},\frac{\pi}{2}]\trans,[\frac{3\pi}{2},\frac{3\pi}{2}]\trans\}$ with $f(\vm x^{*|n})=0$ for $n\in\{1,2,3\}$. Note that the local minima and objective functions satisfy Assumptions \ref{ass:regularity} and~\ref{ass:compatibility}, respectively.
	Figure~\ref{fig:conv_steps} compares SBDP iteration profiles for initializations inside (blue) and outside (red) the theoretical convergence radius from Remark~\ref{remark:estimate_convergenceRadius}, highlighting that the practical convergence region is often larger than the bound of Theorem~\ref{th:linconv}.
	Furthermore, the left part of Figure \ref{fig:conv_curves} shows the convergence to the respective local minima and the right part the empirical convergence rate $C^q:=\frac{\|\vm x\inds - \vm x^{q+1}\|}{\|\vm x\inds - \vm x^{q}\|}$ for the initializations within the theoretical convergence radius, i.e., those marked in blue. Specifically, they are given as $\vm x^{0|n} = \vm x^{*|n}+ [0.25,0.25]\trans$, $n=\{1,2\}$ and $\vm x^{0|3} = \vm x^{*|3}- [0.25,0.25]\trans$.  Remarkably, the rate and order of convergence vary significantly depending on which local minimum SBDP converges to.
	This behavior is explained by the Jacobian \eqref{eq:Jacobian}, with $\vm x = \vm p$
	\begin{equation}
		\vm J(\vm x) = 2\begin{bmatrix}
			\frac{\sin(x_2)}{2 - \sin(x_1)x_2^2}                 & \frac{x_2\cos(x_1)+x_1\cos(x_2)}{2 - \sin(x_1)x_2^2} \\
			\frac{x_2\cos(x_1)+x_1\cos(x_2)}{2 - \sin(x_2)x_1^2} & \frac{\sin(x_1)}{2 - \sin(x_2)x_1^2}
		\end{bmatrix}
	\end{equation}
	for NLP \eqref{eq:example_NLP} which is Lipschitz continuous at each minimum. For $\vm  x^{{*|1}}$, we have $\|\vm J(\vm x^{{*|1}})\| = 0$ such that Theorem \ref{th:linconv} part ii) is applicable near $\vm x^{*|1}$, predicting quadratic convergence. This is confirmed by ${\{C^q\}\rightarrow 0}$ as $\|\vm x^{q,1} - \vm x^{*|1}\| \rightarrow 0$, which is depicted in the right part of Figure~\ref{fig:conv_curves}. Furthermore, the error progressions of $\vm x^{0,2}$ and $\vm x^{0,3}$ display linear convergence. Specifically, we have $\|\vm J(\vm x^{2,*})\|= 0.45 \neq 0$ and $\|\vm J(\vm x^{3,*})\| = 0.08 \neq 0 $, such that Statement i) of Theorem~\ref{th:linconv} applies. In fact, it can be seen that $C^q$ is lower bounded by $\|\vm J(\vm x^{*})\|$ which is a direct consequence of $C\approx \|\vm J(\vm x\inds)\|$ in the neighborhood of $\vm x\inds$, cf. \eqref{eq:ineqJ}. This estimate becomes exact at the optimal solution, explaining the differences in convergence rates $C^q$.
	Figure \ref{fig:compare_curves} shows a comparison of the convergence behavior in vicinity of the local minimum $\vm x_i^{*|1}= \vm 0$ of SBDP with ADMM \cite{Boyd} and ALADIN \cite{Houska} applied to NLP \eqref{eq:example_NLP}  in terms of convergence w.r.t. iterations and communication effort. The initial iterates of all algorithms are initialized at the same distance from the optimum. For ALADIN we use the default implementation offered by the toolbox ALADIN-$\alpha$ \cite{Engelmann3}. The penalty parameter of ADMM is tuned to $1$ for best convergence in the set $\{0.1, 1, 10\}$. Although all three methods solve small-scale constrained NLPs locally, SBDP avoids local copies and has smaller decision vectors. Moreover, SBDP and ADMM require only local computations and neighbor-to-neighbor communication, whereas ALADIN relies on a centralized coordination quadratic program (QP) as well as centralized communication and is hence not fully distributed. The strength of SBDP becomes evident in terms of error progression vs. communicated floats as the agents exchange only one float per iteration (versus two for ADMM and six for ALADIN) for this problem. While ALADIN usually has the best performance since it offers local quadratic convergence guarantees, its performance comes at the high price of centralized computations and communication. Thus, SBDP offers a good trade-off between communication effort and convergence properties compared to state-of-the-art algorithms.
	\subsection{Distributed optimal control of inverted pendulums}
	A frequently used benchmark for the distributed control of nonlinear systems is a chain of $M$ inverted pendulums pivoted to carts which are coupled via springs to their respective neighbors. The nonlinear, continuous-time dynamics of each cart $\agents$ are given by two second-order ODEs
	\begin{align}
		& \ddot \theta_i = \frac{(m_c+m_l)g\sin\theta_i-(m_ll\dot{\theta_i}^2\sin\theta_i + u_i + F_{i})\cos\theta_i}{m_cl + m_ll\sin^2\theta_i}  \nonumber
		\\
		& \ddot y_i = \frac{m_l\sin\theta_i(l \dot{\theta_i}^2  - g \cos\theta_i) + u_i + F_{i} }
		{m_c + m_l \sin^2\theta_i}\label{eq:pendulum_dynamics}
	\end{align}
	with mass $m_c$ of the cart, mass $m_l$ of the load, pendulum length $l$, gravity $g=9.81 \si{\meter \per \second\squared}$, input force $u_i$, displacement $y_i$ and angle $\theta_i$ to the vertical direction. The coupling to the neighboring carts is given by the spring force $F_{i}=F_{i}(y_i,[y_j]_{\neighs}) = \sum_{\neighs} c(y_{j} \scalebox{0.7}[1.0]{\( - \)} y_{i})$ with spring stiffness $c$.
	
In particular, we consider the side-stepping problem over a distance $y_i^T$ in finite time $t\in[0,T]$, leading to the following boundary conditions for all agents: $ y_i(0) = y_i^0,\, y_i(T) = y_i^T,\, \dot{y}_i(0) = \dot{y}_i(T) = 0,\,  \theta_i(0) = 0,\, \theta_i(T) = 0,\, \dot{\theta}_i(0) = \dot{\theta}_i(T) = 0$.
The initial and desired displacements are $y_i^0 = i\scalebox{0.7}[1.0]{\( - \)}1\,\si{\meter}$ and $y_i^T = i\,\si{\meter}$, $\agents$.
Furthermore, we consider the local inequality constraints $y_i\in [y_i^-,\, y_i^+]$, $\dot y_i\in [\dot y^-,\, \dot y^+]$ and $u_i\in \mathbb{U}_i= [u^-,u^+]$ with position bounds $y^{+} =i+1.15\,\si{\meter}$ and $y^{-} = i-0.15 \, \si{\meter}$, velocity limits  $\dot y^{\pm} = \pm 4 \, \si{\meter \per \second}$ and the maximum available force $u^\pm = \pm75\, \si{\newton}$. Let $\mathbb{X}_i$ denote the Cartesian product of these state constraints. We discretize the continuous-time dynamics \eqref{eq:pendulum_dynamics} according to \eqref{eq:DMPC_OCP_equality} using a fourth-order Runge-Kutta scheme to obtain the discretization function $\vm f_i^d(\cdot)$. To preserve the neighbor-affine structure of the dynamics \eqref{eq:pendulum_dynamics} according to \eqref{eq:DMPC_OCP_equality}, not only is the input $u_i^k$ kept constant over the integration interval but also the neighboring states $\vm x_j^k$. We choose $N=80$ intervals leading to an interval length of $\Delta t = 12.5\, \si{\milli\second}$. The discrete-time state and input vectors follow as $\vm x_i^k:=[y_i^k,\dot{y}_i^k,\theta_i^k, \dot{\theta}_i^k]\trans$ and $u_i^k$ where the superscript $k$ denotes the $k$-th discretization point. Finally, we obtain the following NLP
	\begin{subequations}\label{eq:discrete_time_OCP} for all $\agents$ 
		\begin{align}
			\min_{\vm {\hat x },\, \vm{\hat u}} \quad & \sum_{i\in\mathcal V} \sum_{k=0}^{N-1}  (\vm x_i^k)\trans \vm Q \vm x_i^k + R(u_i^k)^2 \label{eq:discrete_time_OCP_costFunction}
			\\  ~\st \quad& \vm x_i^{k+1} =\vm x_i^k\! +\!\Delta t\vm f_i^d(\vm x_i^k,  \vm u_i^k,\Ni{x}^k),\, k \in \mathbb{I}_{[0,N-1]}  \label{eq:discrete_time_OCP_equality} \\
			& \vm x_i^0 = \vm x_0 \,, \smallspace \vm x_i^N = \vm  x_N                                                                                          \\
			& \vm x_i^k \in \mathbb{X}_i\,, \smallspace k \in \mathbb{I}_{[0,N]}\,,  \smallspace u_i^k\in  \mathbb{U}_i, \smallspace k \in \mathbb{I}_{[0,N-1]}
		\end{align}
	\end{subequations}
	with the regularization matrices $\vm Q = \diag(1,1,1,1)$ and ${R=1}$. Notice that \eqref{eq:discrete_time_OCP} is of neighbor-affine structure \eqref{eq:central_NLP_neighboraffine} such that we apply Algorithm \ref{alg:SENSI_with_neighboraffinity}, where we solve the arising local NLPs~\eqref{eq:local_NLP} with sequential quadratic programming (SQP). The local SQP iterations are warm-started with the solution of the previous iteration $q-1$. The problem size is obtained as $n = 4040$, $n_g = 3280$ and $n_h=4860$.
	\begin{figure}[tb]
		\centering
		\includegraphics{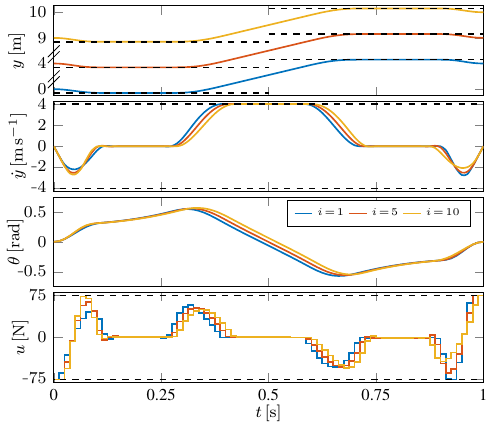}
		\vspace{-4mm}
		\caption{Trajectories of agents $1$, $5$ and $10$ for the side-stepping of multiple coupled inverted pendulums for $ m_c = 1\, \si{\kilogram}$, $m_l=0.25\, \si{\kilogram}$ $l=0.5\, \si{\meter}$ and $c=0.25\, \si{\newton \per \meter}$. The dashed lines indicate the respective constraints.}
		\label{fig:distributed_pendulum}
	\end{figure}
	\begin{figure}[tb]
		\centering
		\includegraphics{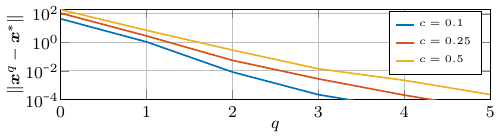}
		\vspace{-4mm}
		\caption{Error of Algorithm \ref{alg:SENSI_with_neighboraffinity} for the example NLP \eqref{eq:discrete_time_OCP} for varying $c \in \{0.1,0.25,0.5\}$.}
		\label{fig:conv_curves_pendulum}
	\end{figure}
	Figure \ref{fig:distributed_pendulum} shows the resulting trajectories for $M=10$ agents with a transition time of $T=1\, \si{\second}$. Algorithm \ref{alg:SENSI_with_neighboraffinity} is initialized with the decentralized solution obtained from a local linear interpolation between the boundary conditions. Note that only the trajectories of the first, middle, and last agent are shown due to the symmetry of the problem. In addition, Figure \ref{fig:conv_curves_pendulum} shows the convergence of the optimizer toward the central solution $\vm x\inds$, obtained by solving the NLP \eqref{eq:discrete_time_OCP} centrally with SQP, for varying spring stiffness strengths $c\in\{0.1,0.25,0.5\}$. Clearly, a higher value of $c$ increases the initial error and reduces the convergence rate, as the norm of the Jacobian $\|\vm J(\vm p)\|$ is proportional to the value of the spring constant $c$. Since $\|\vm J(\vm p\inds)\| \neq 0$, the convergence order is linear, underpinning the result of Theorem \ref{th:linconv}.
	\subsection{Influence of the prediction horizon}
	\begin{figure*}[tb]
		\centering
		\includegraphics{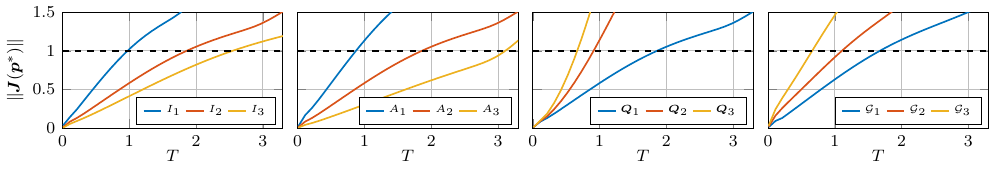}
		\vspace{-4mm}
		\caption{Behavior of the norm of the Jacobian $\| \vm J(\vm p^{*})\|$ w.r.t. to the prediction horizon length $T$ in dependence of different parameter values. In each figure, the value of only one parameter is varied while the other parameters are kept constant at their default value at $I=0.2$, $A=0.1$, $\vm Q=\diag(0,1)$ and $\mathcal{G}_1$. The intersection of the respective curves with the dashed lines, i.e., where $\| \vm J(\vm p^{*})\| = 1$, indicate the maximum allowable prediction horizon length $T_{\mathrm{max}}$.}
		\label{fig:Variation_Jacobian}
	\end{figure*}
	The influence of the prediction horizon length on the convergence of SBDP, when applied in a DMPC context, is evaluated by considering a smart grid consisting of a network of coupled power plants and consumers \cite{Ersdal,Dejong,Burk}. The network is described by two disjunctive sets of generators $\mathcal{V}_G$ and loads $\mathcal{V}_L$ with $\mathcal{V}_G \cup \mathcal{V}_L = \mathcal{V}$. The state is $\vm{ \ubar x}_i = [\varphi\,, \dot \varphi]\trans$, where the phase angle $\varphi(t)$ and angular velocity $\dot \varphi(t)$ are relative to a synchronous equilibrium. The control task consists in stabilizing this equilibrium, i.e., synchronizing all agents at a nominal frequency, e.g., $50\,\si{\hertz}$ in Europe \cite{Ersdal}. The nonlinear, continuous-time dynamics are given in neighbor-affine form
	\begin{equation}\label{eq:smart_grid_dynamics}
	I \ddot{\varphi}_i = - D \dot \varphi_i + u_i + d_i + \sum_{\neighs} A\sin(\varphi_i - \varphi_j) \,,
	\end{equation}
	with moment of inertia $I>0$, friction constant $D>0$ and maximum transferable power $A>0$ (per-unit values). The coupling is given by the power transfer between agent $\agents$ and its neighbors $\neighs$ which depends on the phase shift angle $\varphi_i(t) -\varphi_j(t)$. For generators $i \in \mathcal{V}_G$, the controllable power input is $u_i(t)\in [-0.3,\, 0.3]$ and $d_i=0$, whereas for loads $i \in \mathcal{V}_L$ the input is fixed $u_i(t)=0$ and $d_i(t)=d_i$ denotes the uncontrollable load.
	We discretize the continuous-time dynamics \eqref{eq:smart_grid_dynamics} with a second-order Heun method to obtain the discrete-time dynamics $\vm f_i^d(\cdot)$ in accordance with~\eqref{eq:DMPC_OCP_equality}. We choose the number of discretization intervals as $N = 20$, and design quadratic costs for all $\agents $ as $l_i(\vm{ \ubar x}_i, u_i) = \vm{ \ubar x}_i\trans \vm Q_i \vm{ \ubar x}_i + R_iu_i^2$ and $V_i(\vm{ \ubar x}_i)= \vm{ \ubar x}_i\trans \vm P_i \vm{ \ubar x}_i $ with cost matrices $\vm Q_i=\vm P_i = \diag(0,1)$ and $R_i = 1$ for all $\agents$. Thus, we obtain an NLP of the form \eqref{eq:DMPC_OCP}, where the bounds on the inputs are summarized by the inequality constraints \eqref{eq:DMPC_OCP_inequality}. Furthermore, we consider parameter variations $I_n\in\{0.2,\, 0.4,\, 0.6\}$, $A_n\in \{0.05,0.1,0.2\}  $, $ \vm Q_n \in\{\diag(0,1),\diag(0,2),\diag(0,5)\}$ and $\mathcal{G}_n \in \{\mathcal{G}(\mathcal{V},\mathcal{E}_1), \mathcal{G}(\mathcal{V},\mathcal{E}_2),\mathcal{G}(\mathcal{V},\mathcal{E}_3)\}$, $n \in \{1,2,3\}$, where the values of the varied parameters are similar to the values given in~\cite{Dejong}. The set of edges $\mathcal{E}_n$ describe the connections of a one-hop, two-hop and three-hop ring topology, respectively. We consider an alternating order of generators and loads, i.e., $u_i = 0$ and $d_i = 0.1$, if $i$ is even and $d_i=0$ if $i$ is odd. 
	To illustrate the results of Theorem \ref{th:OCP_case}, we examine the behavior of the norm of the Jacobian in the first DMPC step, $\|\vm J(\vm p\inds)\|$, where $\vm p\inds$ is obtained from a centralized solution with SQP, for a varying prediction horizon $T$. The numerically determined KKT point $\vm p\inds$ satisfies Assumptions~\ref{ass:regularity} and \ref{ass:compatibility} over the considered prediction horizon values $T\in (0,4]$ such that Theorem~\ref{th:OCP_case} is applicable. 
Figure \ref{fig:Variation_Jacobian} shows the norm $\|\vm J(\vm p\inds)\|$ depending on the varying parameters and the prediction horizon. Zero Jacobian terms associated with inactive inequality constraints, i.e., $[\vm \mu\inds]_k = 0$, $k \in \mathbb{I}_{[1,n_h]}$, are omitted. As $T\rightarrow 0$, $\|\vm J(\vm p\inds)\|$ tends toward zero, implying the existence of a maximum prediction horizon for which SBDP converges. This bound depends on the coupling structure in the sense that tightly coupled networks, characterized by a large coupling strength $A$ or many neighbors, reduce the allowable horizon length. Likewise, smaller moments of inertia $I$ or larger $\|\vm Q\|$ decrease $T_{\mathrm{max}}$. In contrast, $\|\vm J(\vm p)\|$ is  invariant w.r.t. the number of agents $M$, provided the coupling topology does not change. Thus, $T_{\mathrm{max}}$ is independent of the network size, which is essential for the scalability of a DMPC scheme.
	\section{Conclusion}
	\label{sec:conclusion}
	This paper introduces a novel SBDP method for the distributed solution of non-convex, large-scale NLPs. It is characterized by solving small-scale NLPs at the agent level and requiring only neighbor-to-neighbor communication. The coordination between agents takes place via sensitivities, which can be computed in a distributed fashion. It is shown that the scheme is at least linearly convergent under appropriate conditions with the potential for quadratic convergence if specific requirements on the coupling structure are met.
	Benefits of the method include applicability to non-convex problems with guaranteed convergence, algorithmic simplicity and a reduced communication effort. However, a limitation is that convergence cannot be guaranteed for all coupling structures. This drawback is mitigated when the SBDP method is considered in a DMPC framework, where it is shown that the prediction horizon can be used to enforce convergence regardless of the problem structure. The theoretical findings are validated in simulations, which also demonstrate favorable convergence properties in practice. Future work investigates modifications to ensure convergence for all coupling structures, different constraint qualifications and the impact of inexact local NLP minimizations.
	%
	%%%%%%%%%%%%%%%%%%%%%%%%%%%%%%%%%%%%%%%%%%%%%%%%%%%%%%%%%%%%%%%%%%%%%%%%%%%%%%%%%%%%%%%%%%%%%%%%%%%%%%%%%%%%%%% Appendix %%%%%%%%%%%%%%%%%%%%%%%%%%%%%%%%%%%%%%%%%%
	%
	\appendix
	\subsection{Explicit representation of the Jacobian} 
	\label{app:Jacobian}
	The single parts of the Jacobian $\vm J(\vm p) = -\vm M^{-1}(\vm p) \vm N(\vm p)$ in \eqref{eq:Jacobian} can be explicitly computed from \eqref{eq:Jacobian_derivation} as the block-diagonal matrix $\vm M(\vm p) = \diag(\vm M_1(\vm p),\dots, \vm M_M(\vm p)) $ with
	\begin{align}\label{eq:M_i}
		\vm M_i(\vm p) = \begin{bmatrix}
			\nablaxx L_i            & \nablax \vm g_i\trans & \nablax \vm h_i\trans \\
			\nablax \vm g_i         & \vm 0                 & \vm 0                 \\
			\vm U_i \nablax \vm h_i & \vm 0                 & \vm h_i
		\end{bmatrix}
	\end{align}
	for each $\agents $ and
	\begin{align} \label{eq:N}
		\vm N(\vm p) =  \begin{bmatrix}
			\vm N_{11}(\vm p)&\dots & \vm N_{1M}(\vm p) \\
			\vdots&\ddots & \vdots \\
			\vm N_{M1}(\vm p)&\dots & \vm N_{MM}(\vm p)
		\end{bmatrix}\,,
	\end{align}
	where the diagonal submatrices in \eqref{eq:N} are defined as
	\begin{align}\label{eq:N_ii}
		\vm N_{ii}(\vm p) = \begin{bmatrix}
			\sum_{\neighs}\nabla_{\vm x_i\vm x_i}^2 L_j            &\vm 0 & \vm 0 \\
			\vm 0         & \vm 0                 & \vm 0                 \\
			\vm 0 & \vm 0                 & \vm 0
		\end{bmatrix}
	\end{align}
	with $i\in \mathcal{V}$ and the off-diagonal matrices as
	\begin{align}\label{eq:N_ij}
		\vm N_{ik}(\vm p) = \begin{bmatrix}
			 \vm L_{ik}  & \nablax \vm g_k\trans & \nablax \vm h_k\trans \\
			\nabla_{\vm x_k} \vm g_i         & \vm 0                 & \vm 0                 \\
			\vm U_i \nabla_{\vm x_k} \vm h_i & \vm 0                 & \vm 0
		\end{bmatrix}
	\end{align}
	with submatrices $\vm L_{ik}:=\nabla_{\vm x_i\vm x_k}^2 L_i + \nabla_{\vm x_k} \sum_{\neighs}\nablax L_j $ and $ k = {1,\dots, M}$.
	\subsection{Proof of Theorem \ref{th:linconv}}
	\label{app:Th1}
	As the Jacobian  is a composition of at least once continuously differentiable functions, there exists some $0<L<\infty$ such that $\|\vm J(\vm p) - \vm J(\vm p\inds)\| \leq L \| \vm p - \vm p\inds\| $ for some $0<r_2$ and $\vm p \in \mathcal{B}_{r_2}(\vm p\inds)$. Let $r_3=\min\{r_1,r_2\}$. Then, for $\vm p\indqp \in \mathcal{B}_{r_3}(\vm p\inds)$, $q=1,2,\dots$, the next SBDP iterate $\vm p\indq$ can be evaluated via the line integral along the linear path $\vm p\inds + s \Delta \vm p\indqp$ with $s\in [0,\,1]$ and $\Delta \vm p\indq := \vm p\indq - \vm p\inds$. The path lies completely within $\mathcal{B}_{r_3}(\vm p\inds)$ by construction and the next iterate $\vm p\indq = \vm \Phi(\vm p\indqp)$ is given as
	\begin{align}\label{eq:path_integralEstimate}
		\vm p\indq = \vm p\inds + \int_0^1 \vm J(\vm p\inds + s\Delta \vm p\indqp)\Delta \vm p\indqp\, \dd s
	\end{align}
	with $\vm p\inds =\vm \Phi(\vm p\inds)$ which follows by the optimality of the distributed solution by Lemma \ref{lem:optConditions} and from $\vm \Phi(\cdot)$ being unique in $\mathcal{B}_{r_3}(\vm p\inds)$ by Lemma \ref{lem:solvability}.
	By adding and subtracting $\vm J(\vm p\inds)\Delta \vm p\indqp$ under the integral on the right-hand-side of \eqref{eq:path_integralEstimate} and taking the Euclidean norm on both sides, we obtain
	\begin{align} \label{eq:ineqJ}
		\| \Delta \vm p\indq\| & \leq \| \vm J(\vm p\inds)\| \| \Delta \vm p\indqp \| \nonumber                                                                         \\
		& \phantom{\leq}+\int_0^1 \vm \|\vm J(\vm p\inds + s\Delta \vm p\indqp) - \vm J(\vm p\inds)\|\|\Delta \vm p\indqp\|\, \dd s \nonumber \\
		& \leq \| \vm J(\vm p\inds)\| \| \Delta \vm p\indqp \| + \frac{L}{2} \|\Delta \vm p\indqp\|^2\,,
	\end{align}
	where the Lipschitz property of $\vm J(\vm q)$ is used in the second inequality. The inequality \eqref{eq:ineqJ} is utilized to prove the individual statements i) and ii). Statement i) can be proven by induction. Let $r$ be such that $r<\min\{r_3, 2(C - \|\vm J(\vm p^*)\|)/L$ with $0<\|\vm J(\vm p^*)\|<C<1$. For $q=1$ (induction start), the relation is satisfied since$\| \vm p^1 - \vm p\inds \| \leq C\| \vm p^0 - \vm p\inds \|$ which follows from the assumption that  $\vm p^0 \in \mathcal{B}_{r}(\vm p\inds)$.
	Postulating that $\vm p^{q-1} \in \mathcal{B}_{r}(\vm p\inds)$ holds for step $q-1$ (induction hypothesis), we get from \eqref{eq:ineqJ} 
	\begin{equation}\label{eq:lin_convappendix}
		\| \Delta \vm p\indq\| \leq \| \vm J(\vm p\inds)\| \| \Delta \vm p\indqp \| +  \frac{L}{2} \|\Delta \vm p\indqp\|^2 \!\leq\!  C\| \Delta \vm p\indq \|
	\end{equation}
	which completes the proof by induction. 
	Thus, if $\vm p^0 \in \mathcal{B}_r(\vm p\inds)$, $r< \min\{r_3, 2(C - \|\vm J(\vm p^*)\|)/L \}$ and $\|\vm J(\vm p^*)\|<1$, then the sequence $\{\vm p\indq\}$ generated by the SBDP method is well-defined as it completely lies within $\mathcal{B}_r(\vm p\inds)$ and converges linearly to a (local) minimizer $\vm p\inds$.
	We continue with Statement ii). If $\|\vm J(\vm p\inds)\|=0$, then from \eqref{eq:ineqJ} we have
	\begin{align} \label{eq:ineqJ2}
			\| \Delta \vm p\indq\| \leq  \frac{L}{2} \|\Delta \vm p\indqp\|^2\,.
	\end{align}
	Thus, Statement ii), i.e., quadratic convergence directly follows from \eqref{eq:ineqJ2} with $\tilde C = L/2$ and the same inductive reasoning as in Statement i). This proves Theorem \ref{th:linconv}.  \hfill \QEDclosed
	\subsection{Proof of Theorem \ref{th:OCP_case}}
	\label{app:Th2}
	We follow along the lines of Theorem~\ref{th:linconv} and inspect the behavior of the Jacobian $\vm J(\vm p\inds)=-\vm M(\vm p\inds)^{-1}\vm N(\vm p\inds)$ w.r.t. the prediction horizon $T$ for a given $N$, resulting in a Jacobian of constant dimension. If we insert $T /N = \Delta t $ into \eqref{eq:DMPC_OCP_costFunction} and \eqref{eq:DMPC_OCP_equality}, we can interpret $T$ as parameter in \eqref{eq:DMPC_OCP} for a constant $N$. Since Assumptions \ref{ass:regularity} and \ref{ass:compatibility} hold for any $ T \in \mathbb{T}$, there exists a well-defined KKT-point $\vm p\inds = \vm p\inds(T)$ of OCP~\eqref{eq:DMPC_OCP}, which now depends on $T$. Thus, for any $\vm p\in \mathcal{B}_r(\vm p\inds(T))$ with $T \in  \mathbb{T}$ the statements of Lemma~\ref{lem:solvability} are applicable and in particular $\vm M( \vm p\inds(T); T) $ is regular. Since $\vm p\inds(T)$ is uniquely determined by $T$, we obtain $\vm M(T)= \vm M( \vm p\inds(T); T)$ and $\vm N(T)= \vm N(\vm p\inds(T); T)$ which only depend on $T$. Due to the structure of the discretized dynamics~\eqref{eq:DMPC_OCP_equality} and costs \eqref{eq:DMPC_OCP_costFunction}, we can partition $\vm M( T)$ as the following expression
	\begin{align}\label{eq:expansion_M}
		\vm M( T) = \vm M_0 + T \vm M_1(T)\,,
	\end{align}
	where $\vm M_0\in \mathbb{R}^{p\times p}$ and $\vm M_1(T) = \vm M_1(\vm p\inds(T)) \in \mathbb{R}^{p\times p} $ are in general singular matrices arising from the specific structure of $\vm M(T)$. However, $\vm M(T)= \vm M_0 + T \vm M_1(T)$ is regular for all $T \in  \mathbb{T}$ according to Lemma~\ref{lem:solvability}. In particular, $\vm M_0$ captures the terms in \eqref{eq:DMPC_OCP} which are independent of $T$, while $\vm M_1(T)$ reflects the coefficients of the linear $T$-terms. The coupling between agents is only present in the running costs \eqref{eq:DMPC_OCP_costFunction} and dynamics \eqref{eq:DMPC_OCP_equality}, such that $ \vm N(T)$ is written as
	\begin{align} \label{eq:expansion_N}
	\vm N(T) = T \vm N_1(T)\,,
	\end{align}
	where $\vm N_1(T) = \vm N_1(\vm p\inds(T))\in \mathbb{R}^{p\times p}$ implicitly depends on $T$ via $\vm p\inds(T)$. The partitionings \eqref{eq:expansion_M} and \eqref{eq:expansion_N} rely on the fact the $T$ explicitly appears via $\Delta t$ and $N$ in the discretized OCP \eqref{eq:DMPC_OCP}. 
	In the following, we establish a condition on $T$ such that 
	\begin{align}\label{eq:bound_jacobian}
		\|\vm J(T )\|= \|(\vm M_0 + T \vm M_1(T))^{-1}T \vm N_1(T)\| < 1
	\end{align}
	with $\vm J(T) = \vm J(\vm p\inds(T);T)$ which is sufficient for convergence according to Theorem \ref{th:linconv}.
	To this end, consider the singular value decomposition (SVD) of $\vm M_0$, i.e., $\vm M_0 = \vm U_0 \vm \Sigma_0 \vm V_0\trans$, where  $\vm \Sigma_0 = \diag(\sigma_{1},\dots, \sigma_{r}, 0,\dots, 0)$ is a diagonal matrix of positive singular values with $\sigma_{i}>0$ for $i=1,\dots,r$, $\sigma_{i}=0$ for $i=r+1,\dots,p$ with $r:=\rank(\vm M_0)$. We now interpret $T \vm M_1(T)$ as a perturbation of $\vm M_0$ for small $T>0$. The perturbed singular values $\sigma_i(T)>0$ of $\vm M_0$ can be computed as the following expansion
	\begin{subequations}
	\begin{align}\label{eq:linear_approx_singular_values}
		\sigma_i(T) &= \sigma_{i} + \Delta \sigma_i T + k_i(T )\,, && i=1,\dots,r \\
		\sigma_i(T) &= \Delta \sigma_i T + k_i(T)\,, && i=r+1,\dots,p\label{eq:linear_approx_singular_values2}
	\end{align}
	\end{subequations}
which is a well known result from perturbation theory \cite{Soderstrom,Stewart}. 
The explicit computation of the linear perturbation values $\Delta \sigma_i $ can be found in \cite[Thm. 4.1 and 4.2]{Soderstrom} and the functions $ k_i(\cdot)\in \mathcal{O}(T^2)$ summarize the higher-order terms. Furthermore, let $\nullspace(\vm M_0):= \{\vm y \in \mathbb{R}^p| \vm M_0\vm y = \vm 0 \} $ and $\nullspace(\vm M_1):= \{\vm y \in \mathbb{R}^p| \vm M_1\vm y = \vm 0 \} $ define the respective null spaces. An important consequence of the regularity of $\vm M$ is that $\nullspace(\vm M_0) \cap \nullspace(\vm M_1)= \{\,\}$ which implies that $\Delta \sigma_i> 0$ for $\sigma_{i}= 0$, $i=r+1,\dots,p$. With the SVD of $\vm M_0$ and the perturbed singular values \eqref{eq:linear_approx_singular_values} and \eqref{eq:linear_approx_singular_values2}, we rewrite \eqref{eq:bound_jacobian} as
	\begin{align}
		\|\vm J(T)\| &= \bigg\|\vm V_0 \begin{bmatrix}
			\vm \Sigma_{r}(T) & \vm 0 \\
			\vm 0 & \vm \Sigma_{p-r}(T)
		\end{bmatrix}^{-1} \vm U_0\trans T \vm N_1(T)\bigg \|\label{eq:bound_jacobian1} \\ &= \bigg\|\vm V_0  \begin{bmatrix}
			T\vm \Sigma_{r}^{-1}(T) & \vm 0 \\
			\vm 0 & T \vm \Sigma_{p-r}^{-1}(T)
		\end{bmatrix} \begin{bmatrix}
	    \vm N_{1,r}(T) \\ \vm 0 
		\end{bmatrix} \bigg\| \nonumber 
	\end{align}
	with the matrices $ \vm \Sigma_r(T) := \diag(\sigma_1(T), \dots,\sigma_r(T) )$ and $ \vm \Sigma_{p-r}(T) := \diag(\sigma_{r+1}(T), \dots,\sigma_p(T) )$ corresponding to the singular values of the range and null space, respectively. Furthermore, $\vm N_{1,r}(T)$ denotes the projection of $\vm N_{1}(T)$ onto the range of $\vm M_0$, i.e., $\vm U_0\trans \vm N_{1}(T) = [\vm N_{1,r}(T)\trans\,, \vm 0\trans ]\trans$.  An upper bound to \eqref{eq:bound_jacobian1} can be found by lower bounding the respective singular values $\sigma_i(T)$ from \eqref{eq:linear_approx_singular_values}. The notation $\mathcal{O}(T^2)$ implies that there exist constants $\delta_1,c_i>0$ such that for any $k_i(\cdot) \in \mathcal{O}(T^2)$ it follows that $|k_i(T)| \leq c_iT^2$ for  $T\in(0,\delta_1]$. Let $\delta_2<\min\{\delta_1, |\Delta \sigma_1|/ c_1,\dots,|\Delta \sigma_p|/ c_p \}$, then the perturbation expansion~\eqref{eq:linear_approx_singular_values} of the first $i=1,\dots,r$ non-zero singular values is lower bounded for $T \in (0,\delta_2]$ as
	\begin{align} \label{eq:lower_bound_nonzerosingularvalues}
		\sigma_i(T) &\geq \sigma_{i}  -|\Delta \sigma_i T|  -|k_i(T)| \geq \sigma_{i} - |\Delta \sigma_i| T - c_iT^2 \nonumber
		\\&=\sigma_{i} - T(|\Delta \sigma_i| + c_iT)\geq \sigma_{i} - 2|\Delta \sigma_i| T   \\
		& \geq \min_{i \in \{1,\dots,r\}} \{\sigma_{i} - 2|\Delta \sigma_i| T\} \geq \sigma_{\min} - 2\Delta \sigma_{\mathrm{max}} T\nonumber
	\end{align}
with the abbreviations $\sigma_{\min} := \min_{i\in\{1,\dots,r\}} \sigma_i$ and $\Delta \sigma_{\mathrm{max}} := \max_{i\in\{1,\dots,r\}} |\Delta \sigma_i|$. Similarly, the perturbation expansion \eqref{eq:linear_approx_singular_values2} of the last $i=r+1,\dots,p$ zero singular values are bounded for $T \in (0,\delta_2]$ as
	\begin{align}\label{eq:lower_bound_zerosingularvalues}
		\sigma_i(T)  &\geq \Delta \sigma_i T - c_iT^2 \geq \Delta \sigma_i (1-T)T\\
		&\geq \min_{i \in \{r+1,\dots,p\}}\!\{ \Delta \sigma_i (1-T)T\}  = \Delta \sigma_{\mathrm{min}} (1-T)T \nonumber 
	\end{align}
with $\sigma_{\mathrm{min}} := \min_{i\in\{r+1,\dots,p\}} |\Delta \sigma_i|$ and where in the first inequality it is used that $\Delta \sigma_i$ is necessarily positive due to the regularity of $\vm M$. For large $T$ both lower bounds get significantly less accurate. Thus, it holds for all $T \in (0,\delta_2]$
	\begin{subequations} \label{eq:upper_bound_diagonal_singular_values}
	\begin{align}
	\|T \vm \Sigma_{r}^{-1}(T)\|	& \leq \frac{T}{\sigma_{\min} - 2\Delta \sigma_{\mathrm{max}}T} \\
		\|T \vm \Sigma_{p-r}^{-1}(T)\| & \leq \frac{1}{\Delta \sigma_{\mathrm{min}} (1-T)}\,. \label{eq:upper_bound_diagonal_singular_values2}
	\end{align}
		\end{subequations}
This essentially shows that the singular values corresponding to the range of $\vm M$ tend to zero for $T \rightarrow 0$, while the singular values corresponding to the null space remain bounded and finite. 
	Inserting the estimates in \eqref{eq:upper_bound_diagonal_singular_values} into \eqref{eq:bound_jacobian1} leads to 
	\begin{align}
		\|\vm J(T)\| &\leq \bigg\|\vm V_0  \begin{bmatrix}
			T\vm \Sigma_{r}^{-1}(T) & \vm 0 \\
			\vm 0 & T \vm \Sigma_{p-r}^{-1}(T)
		\end{bmatrix} \begin{bmatrix}
				 \vm {N}_{1,r}(T) \\ \vm 0 
		\end{bmatrix} \bigg \| \nonumber  \\
		&\leq \frac{T \|\vm V_{0}\|\|	 \vm {N}_{1,r}(T)\|}{\sigma_{\min} - 2\Delta \sigma_{\mathrm{max}}T} \leq \frac{T  \|\vm {N}_{1,r}(T)\|}{\sigma_{\min} - 2\Delta \sigma_{\mathrm{max}}T}\,, 
	\end{align}
	where  $\|\vm V_{0}\|=1$ is used. Note that $\vm {N}_{1,r}(T) = \vm {N}_{1,r}(\vm p\inds(T))$ is still implicitly dependent on $T$. However, we can upper bound $\|\vm {N}_{1,r}(\vm p\inds(T))\|$ by considering the maximum value $n_{\mathrm{max}}:=\max_{T \in \mathbb{T}} \|\vm {N}_{1,r}(\vm p\inds(T))\|$ to obtain 
	\begin{align}
	\|\vm J(T)\|\leq \frac{T n_{\mathrm{max}} }{\sigma_{\min} - 2\Delta \sigma_{\mathrm{max}}T}
	\end{align}
    which is guaranteed to exist since the interval $\mathbb{T}$ is compact and $\|\vm N_{1,r}(\vm p\inds(T))\|$ is continuous.  
	Imposing the condition $\|\vm J(T)\|< 1$ for convergence according to Theorem \ref{th:linconv}, results in an upper bound on the prediction horizon $T$
	\begin{align} \label{eq:upper_bound_delta_t}
	T< \frac{\sigma_{\min}}{ n_{\mathrm{max}} +2\Delta \sigma_{\mathrm{max}}}=:T_{\mathrm{max}}\,.
	\end{align}
	Note that $ \sigma_{\mathrm{min}} \neq 0$ holds since $\vm M$ is regular establishing that $T_{\mathrm{max}}>0$. Furthermore, we need to impose that $(0,\delta_2]\cap\mathbb{T} \neq \{\,\}$ which guarantees that the bound~\eqref{eq:upper_bound_delta_t} is attainable. Finally, by imposing that $T<T_{\mathrm{max}}$ and $T_{\mathrm{max}} \in \mathrm{int}(\mathbb{T})$, ensures 
	convergence of SBDP by Theorem~\ref{th:linconv} proving Theorem~\ref{th:OCP_case}. \hfill \QEDclosed 
	\vspace{-3mm}
	%%%%%%%%%%%%%%%%%%%%%%%%%%%%%%%%%%%%%%%%%%%%%%%%%%%%%%%%%%%%%%%%%%%%%%%%%%%%%%%%%%%%%%%%%%%%%%%%%%%%%%%%%%%%%%%%
	\bibliographystyle{IEEEtran}
	\bibliography{Pierer_TCNS_bib_2026}
	%%%%%%%%%%%%%%%%%%%%%%%%%%%%%%%%%%%%%%%%%%%%%%%%%%%%%%%%%%%%%%%%%%%%%%%%%%%%%%%%%%%%%%%%%%%%%%%%%%%%%%%%%%%%%%%%
\end{document}